\documentclass{gtart_h}

\def\ifplaintex{\expandafter\ifx\csname documentclass\endcsname\relax}

\def\gtp{{\mathsurround=0pt\it $\cal G\mskip-2mu$eometry \&\ 
$\cal T\!\!$opology $\cal P\!$ublications}}  

\def\recd{{\small Received:\qua\receiveddate\ifx\reviseddate\relax
\else\qquad Revised:\qua\reviseddate\fi\par}} 

\def\lognumber#1{\def\thelognumber{#1}}
\def\volumenumber#1{\def\thevolumenumber{#1}}
\def\volumeyear#1{\def\thevolumeyear{#1}}
\def\papernumber#1{\def\thepapernumber{#1}}
\def\pagenumbers#1#2{\def\startpage{#1}\def\finishpage{#2}}
\def\published#1{\def\publishdate{#1}}

\def\received#1{\def\receiveddate{#1}}

\def\accepted#1{\def\accepteddate{#1}}

\def\asciiaddress#1{\def\theasciiaddress{#1}}
\def\asciiemail#1{\def\theasciiemail{#1}}

\long\def\asciiabstract#1{\long\def\theasciiabstract{#1}}

\let\\\par\let\thelognumber\relax\let\thevolumenumber\relax
\let\thepapernumber\relax\let\thevolumeyear\relax\let\startpage\relax
\let\finishpage\relax\let\publishdate\relax\let\receiveddate\relax
\let\reviseddate\relax\let\accepteddate\relax\let\theasciititle\relax
\let\theasciiauthors\relax\let\theasciiaddress\relax
\let\theasciiabstract\relax

\let\theasciiemail\relax

\ifplaintex
\font\logobig=cmssbx10 scaled 3836
\font\logomed=cmssbx10 scaled 2557
\else
\font\logobig=cmssbx10 scaled 4200
\font\logomed=cmssbx10 scaled 2800
\fi

\long\def\makeagttitle{   
\count0=\startpage
\agt\hfill      
\hbox to 45truept{\vbox to 0pt{\vglue -13truept{\logomed A\kern -.37em{\logobig 
T}\kern -.38em G}\vss}\hss}
\break
{\small Volume \thevolumenumber\ (\thevolumeyear)
\startpage--\finishpage\nl
Published: \publishdate}

\vglue .25truein

{\parskip=0pt\leftskip 0pt plus
1fil\def\\{\par\smallskip}{\Large\bf\thetitle}\par\medskip} \vglue
0.05truein

{\parskip=0pt\leftskip 0pt plus 1fil\def\\{\par}{\sc\theauthors}
\par\medskip}%
 
\vglue 0.03truein 

{\small\leftskip 25truept\rightskip 25truept{\bf Abstract}\stdspace\theabstract

{\bf AMS Classification}\stdspace\theprimaryclass
\ifx\thesecondaryclass\relax\else; \thesecondaryclass\fi\par
{\bf Keywords}\stdspace \thekeywords\par}\vglue 7truept

}   

\ifplaintex
\font\phead=cmsl9 scaled 950
\font\pnum=cmbx10 scaled 913
\font\pfoot=cmsl9 scaled 950
\headline{\vbox to 0pt{\vskip -4.5mm\line{\small\phead\ifnum
\count0=\startpage ISSN 1472-2739 (on-line) 1472-2747 (printed)
\hfill {\pnum\folio}\else\ifodd\count0\def\\{ }%
\ifx\theshorttitle\relax\thetitle\else\theshorttitle\fi\hfill{\pnum\folio}
\else\def\\{ and }{\pnum\folio}\hfill\ifx\theshortauthors\relax\theauthors
\else\theshortauthors\fi\fi\fi}\vss}}
\footline{\vbox to 0pt{\vglue 0mm\line{\small\pfoot\ifnum\count0=\startpage
\copyright\ \gtp\hfill\else
\agt, Volume \thevolumenumber\ (\thevolumeyear)\hfill\fi}\vss}}
\else
\font=cmsl9 scaled 1050
\font\lnum=cmbx10 
\font=cmsl9 scaled 1050
\makeatletter
\def\@oddhead{{\small\lhead\ifnum\count0=\startpage ISSN 1472-2739 
(on-line) 1472-2747 (printed)\hfill {\lnum\number\count0}\else\ifodd\count0
\def\\{ }\ifx\theshorttitle\relax \thetitle \else\theshorttitle\fi\hfill
{\lnum\number\count0}\else\def\\{ and }{\lnum\number\count0}
\hfill\ifx\theshortauthors\relax 
\theauthors\else\theshortauthors\fi\fi\fi}}\def\@evenhead{\@oddhead}
\def\@oddfoot{\small\lfoot\ifnum\count0=\startpage\copyright\ \gtp\hfill\else
\agt, Volume \thevolumenumber\ (\thevolumeyear)\hfill\fi}
\def\@evenfoot{\@oddfoot}
\makeatother
\fi
\let\maketitlepage\makeagttitle
\let\makeshorttitle\maketitlepage
\let\maketitle\maketitlepage

\newwrite\gtoutfile
\long\gdef\makeheadfile{  
{\def\\{, }\def\s{ }
\immediate\openout\gtoutfile head.xxx
\immediate\write\gtoutfile{Proxy-for: \ifx\theasciiauthors\relax
\theauthors\else\theasciiauthors\fi\s<\ifx\theasciiemail\relax\theemail\else\theasciiemail\fi>}
\immediate\write\gtoutfile{\noexpand\\}
\immediate\write\gtoutfile{Authors: \ifx\theasciiauthors\relax
\theauthors\else\theasciiauthors\fi}
{\def\\{ }\immediate\write\gtoutfile{Title: \ifx\theasciititle\relax
\thetitle\else\theasciititle\fi}}
\immediate\write\gtoutfile{Subj-class: GT or SG, GR etc}
\immediate\write\gtoutfile{MSC-class: \theprimaryclass\ifx\thesecondaryclass\relax\else, \thesecondaryclass\fi}
\immediate\write\gtoutfile{Journal-ref: Algebr. Geom. Topol. \thevolumenumber\s
(\thevolumeyear) \startpage-\finishpage}
\immediate\write\gtoutfile{Comments: Published by Algebraic and
Geometric Topology at}
\immediate\write\gtoutfile{\s\s\s  http://www.maths.warwick.ac.uk/agt/AGTVol\thevolumenumber/agt-\thevolumenumber-\thepapernumber.abs.html}
\immediate\write\gtoutfile{\noexpand\\}
\immediate\write\gtoutfile{}
\ifx\theasciiabstract\relax
\immediate\write\gtoutfile{\theabstract}\else
\immediate\write\gtoutfile{\theasciiabstract}\fi
\immediate\write\gtoutfile{}
\immediate\write\gtoutfile{\noexpand\\}
\immediate\write\gtoutfile{}
\immediate\closeout\gtoutfile}}  

\def\maketitlepage{\makeagttitle\makeheadfile}
\let\makeshorttitle\maketitlepage
\let\maketitle\maketitlepage

\lognumber{14}
\volumenumber{5}
\volumeyear{2005}
\papernumber{14}
\pagenumbers{301}{354}
\received{9 October 2004} 
\accepted{8 March 2005}
\published{17 April 2005}

\usepackage{amssymb,amsmath,graphicx,psfrag}
\def\psfraga <#1,#2> #3#4{%
\psfrag {#3}{\smash{\rlap{\kern #1 \raise #2\hbox{#4}}}}}

\def\S{Section }

\newcommand{\mc}[1]{{\mathcal #1}}

\newtheorem{theorem}{Theorem}[section]
\newtheorem{proposition}[theorem]{Proposition}
\newtheorem{corollary}[theorem]{Corollary}
\newtheorem{lemma}[theorem]{Lemma}

\theoremstyle{definition}
\newtheorem{definition}[theorem]{Definition}
\newtheorem{remark}[theorem]{Remark}

\newcommand{\pin}{p_{\op{in}}}
\newcommand{\pout}{p_{\op{out}}}

\newcommand{\eqdef}{\;{:=}\;}
\newcommand{\fedqe}{\;{=:}\;}

\newcommand{\Q}{{\mathbb Q}}
\newcommand{\R}{{\mathbb R}}
\newcommand{\Z}{{\mathbb Z}}
\newcommand{\op}{\operatorname}

\newcommand{\M}{\mc{M}}

\newcommand{\Ker}{\op{Ker}}

\newcommand{\tensor}{\otimes}

\newcommand{\feed}{-\!\!\!\hspace{0.7pt}\mbox{\raisebox{3pt}{$\shortmid$}}\,}

\newcommand{\union}{\bigcup}

\begin{document}

\title{The periodic Floer homology of a Dehn twist}
\authors{Michael Hutchings\\Michael Sullivan}
\address{Department of Mathematics,
University of California\\Berkeley, CA 94720-3840, USA }
\secondaddress{Department of
Mathematics and Statistics, University of Massachusetts\\Amherst,
MA 01003-9305, USA}
\asciiaddress{Department of Mathematics,
University of California\\Berkeley CA 94720-3840, USA\\and\\Department of
Mathematics and Statistics, University of Massachusetts\\Amherst,
MA 01003-9305, USA}
\asciiemail{hutching@math.berkeley.edu, sullivan@math.umass.edu}
\gtemail{\mailto{hutching@math.berkeley.edu}{\rm\qua and\qua}\mailto{sullivan@math.umass.edu}}

\begin{abstract}
  The periodic Floer homology of a surface symplectomorphism, defined
  by the first author and M. Thaddeus, is the homology of a chain
  complex which is generated by certain unions of periodic orbits, and
  whose differential counts certain embedded pseudoholomorphic curves
  in $\R$ cross the mapping torus.  It is conjectured to recover the
  Seiberg-Witten Floer homology of the mapping torus for most spin-c
  structures, and is related to a variant of contact homology.  In
  this paper we compute the periodic Floer homology of some Dehn
  twists.
\end{abstract}
\asciiabstract{%
  The periodic Floer homology of a surface symplectomorphism, defined
  by the first author and M. Thaddeus, is the homology of a chain
  complex which is generated by certain unions of periodic orbits, and
  whose differential counts certain embedded pseudoholomorphic curves
  in R cross the mapping torus.  It is conjectured to recover the
  Seiberg-Witten Floer homology of the mapping torus for most spin-c
  structures, and is related to a variant of contact homology.  In
  this paper we compute the periodic Floer homology of some Dehn
  twists.}

\primaryclass{57R58}\secondaryclass{53D40, 57R50}
\keywords{Periodic Floer homology, Dehn twist, surface symplectomorphism}

\maketitle

\section{Introduction}
\label{sec:introduction}

Let $\Sigma$ be a compact surface, possibly with boundary, with a
symplectic form $\omega$, and let $\phi$ be a symplectomorphism of
$\Sigma$.  Consider the {\em mapping torus\/}
\begin{equation}
\label{eqn:mappingTorus}
Y_\phi\eqdef\frac{[0,1]\times \Sigma}{(1,x)\sim (0,\phi(x))}.
\end{equation}
Projection onto the $[0,1]$ factor defines a fibration $Y_\phi\to
\R/\Z$ with fiber $\Sigma$.  There is a natural flow $R$ on $Y_\phi$
which increases the $[0,1]$ coordinate, and we identify periodic
orbits of $\phi$ with closed orbits of this flow in $Y_\phi$.  We fix
a homology class $h\in H_1(Y_\phi)$, and define the {\em degree\/} $d$
to be the intersection number of $h$ with a fiber of $Y_\phi\to\R/\Z$.
Under ``monotonicity'' and ``$d$-regularity'' assumptions on $\phi$,
we can define the {\em periodic Floer homology\/} (PFH), denoted by
$HP_*(\phi,h)$; see \cite{pfh1}, the background in \cite{pfh2}, and
the review in \S\ref{sec:pfh}.  Roughly speaking, this is the homology
of a chain complex which is generated by certain unions of periodic
orbits of $\phi$ with total homology class $h$, and whose differential
counts certain embedded pseudoholomorphic curves in $\R\times
Y_{\phi}$ for a suitable almost complex structure $J$.  It is expected
that the PFH is independent of $J$ and invariant under appropriate
isotopy of $\phi$ fixing $\phi|_{\partial\Sigma}$.

When $d=1$, PFH reverts to the ordinary symplectic Floer homology of
$\phi$, which is the homology of a chain complex generated by fixed
points of $\phi$ and whose differential counts pseudoholomorphic
sections of $\R\times Y_{\phi}\to\R\times S^1$.  The symplectic Floer
homology has been computed for a Dehn twist by Seidel \cite{s96}, more
examples have been computed by Eftekhary \cite{ee} and Gautschi
\cite{rg}, and a conjecture for a product of positive Dehn twists on a
surface with boundary is given in \cite{s00}.

\paragraph{Results}
In this paper we compute the periodic Floer homology of some Dehn
twists.  We begin, in Theorem~\ref{thm:cylinder}, by computing the PFH
of a perturbation of the ``positive'' twist on the cylinder
\begin{equation}
\label{eqn:cylinderTwist}
\begin{split}
\phi_0:[X_1,X_2]\times S^1 & \longrightarrow [X_1,X_2]\times S^1,\\
(x,y) & \longmapsto (x,y-x)
\end{split}
\end{equation}
for any $J$.  Here, as elsewhere in this paper, we identify
$S^1=\R/\Z$.  In fact, in Theorem~\ref{thm:rounding} we obtain a
combinatorial formula for most of the differential in the chain
complex, in terms of ``rounding corners'' of convex polygonal paths
connecting lattice points in the plane.  In the terminology of
\S\ref{sec:pfh}, this part of the differential counts genus zero
pseudoholomorphic curves with two ``incoming ends'', any number of
``outgoing ends'', possibly together with ``trivial cylinders''.  We
prove these theorems in \S\ref{sec:cylinder}.

In Theorem~\ref{thm:torus}, proved in \S\ref{sec:torus}, we compute
the PFH of an $n$-fold positive Dehn twist on the torus (for any $J$),
namely a perturbation of
\begin{equation}
\label{eqn:torusTwist}
\begin{split}
\phi_0^T: (\R/n\Z)\times (\R/\Z) & \longrightarrow (\R/n\Z)\times (\R/\Z),\\
(x,y) & \longmapsto (x,y-x).
\end{split}
\end{equation}
In \S\ref{sec:surface}, we consider a composition $\phi^\Sigma$ of
Dehn twists along disjoint circles $\gamma_i$ on a higher genus
(connected) surface $\Sigma$, perturbed so that it is close to the
identity away from the circles $\gamma_i$.  To ensure monotonicity of
the standard representative of the composition of Dehn twists, we
assume:
\begin{itemize}
\item[$(*)$] If $\xi\in H_1(\Sigma)$ has nonzero intersection number
  with some $[\gamma_i]$, then $\phi^\Sigma_*\xi\neq\xi$ in
  $H_1(\Sigma)$.  Also, $d\neq g(\Sigma)-1$ when
  $\partial\Sigma=\emptyset$.
\end{itemize}
If we further assume that
\begin{itemize}
\item[$(**)$] each component of $\Sigma\setminus\bigcup_i\gamma_i$ contains
  a component of $\partial\Sigma$ or has sufficiently large genus with
  respect to $d$,
\end{itemize}
then the PFH can be computed in terms of the cylinder complex
described above and Morse theory on $\Sigma\setminus\bigcup_i\gamma_i$.
This is complicated in general, but we work out the cases when
$\Sigma$ is closed and there is one circle $\gamma_i$, nonseparating
or separating, in Theorems~\ref{thm:nonseparating} and
\ref{thm:separating} respectively.

\paragraph{Motivation}
There are three basic motivations for this paper.  First, it is
conjectured in \cite{pfh1} that if $\Sigma$ is closed and connected
and $d<g(\Sigma)-1$, then $HP_*(\phi,h)$ is isomorphic to the
Seiberg-Witten Floer (SWF) homology of $Y_\phi$ for a spin-c structure
corresponding to $h$.  This conjecture is an analogue of Taubes's
``SW=Gr'' theorem \cite{t95}, relating Seiberg-Witten solutions to
embedded pseudoholomorphic curves on closed symplectic 4-manifolds,
for the noncompact symplectic 4-manifold $\R\times Y_\phi$.  This
conjecture is known to hold at the level of Euler characteristics as a
consequence of Taubes's theorem applied to $S^1\times Y_\phi$, see
\cite{hl1}.  It also fits nicely with a conjecture of Salamon
\cite{sal} relating SWF homology of mapping tori to symplectic Floer
homology of induced maps on symmetric products.  Thus our results
conjecturally give some new examples of SWF homology.  The version of
SWF homology considered here is conjectured to agree with the
invariant $HF^+$ of Ozsv\'{a}th-Szab\'{o} \cite{os1,os2}.

Second, one should be able to define a ``quantum product'' relating
the PFH of three surface diffeomorphisms $f$, $g$, and $f\circ g$,
cf.\ \cite{skd}, and then in principle recover the Seiberg-Witten
invariants of a 4-dimensional symplectic Lefschetz pencil by taking
the quantum products of certain elements in the PFH of a Dehn twist.
Thus the PFH of a Dehn twist appears to be a fundamental building
block in invariants of symplectic 4-manifolds.

Third, one can define an analogue of PFH for a contact 3-manifold,
thus obtaining an interesting variant of the symplectic field theory of
\cite{egh}, which we call ``embedded contact homology''.  Some basic
contact 3-manifolds, such as $S^1\times S^2$ with the overtwisted
contact structure considered in \cite{t01}, or the unit cotangent
bundle of $T^2$, contain pieces whose Reeb flow is diffeomorphic to the
mapping torus flow $R$ for the inverse of the cylinder twist
\eqref{eqn:cylinderTwist}.  Thus we expect the techniques developed in
this paper to be useful in computing more examples of various flavors
of contact homology.

\paragraph{Remarks on the proofs}
We organize the computations in this paper using some spectral
sequences, most of which are variations on the following theme.
Suppose that our symplectomorphism $\phi$ has an invariant curve
$\xi\subset\Sigma$, on which $\phi$ is an irrational rotation, and
which divides $\Sigma$ into pieces $\Sigma_1$ and $\Sigma_2$.  Then in
some cases, there is a spectral sequence $\mc{E}^*_{*,*}$ which
converges to $HP_*(\phi,h)$, with
\begin{equation}
\label{eqn:invariantCurve}
\mc{E}^1\simeq \bigoplus_{h_1+h_2=h}
HP_*\left(\phi|_{\Sigma_1},h_1\right)
\tensor
HP_*\left(\phi|_{\Sigma_2}, h_2\right).
\end{equation}
One such spectral sequence is used in \S\ref{sec:cylinderSpectral} to
compute the PFH of a twist on a cylinder, by studying how the PFH
changes as one twists the symplectomorphism on the boundary.  Other
such spectral sequences are used for the torus in \S\ref{sec:eta} and
for higher genus surfaces in \S\ref{sec:nonseparating} and
\S\ref{sec:separating} to cut out the cylinders in which twisting
takes place.

In some cases a spectral sequence satisfying
\eqref{eqn:invariantCurve} exists but might not converge to
$HP_*(\phi,h)$, roughly speaking due to higher genus pseudoholomorphic
curves which ``wrap around'' large parts of $\Sigma$.  This failure to
converge is measured by the ``wrapping spectral sequence'' introduced
in \S\ref{sec:eta}.  We make the hypothesis $(**)$ in order to avoid
dealing with these higher genus curves by ensuring that the wrapping
spectral sequence degenerates.

An simple but important trick used in several places is the ``local
energy inequality'' of \S\ref{sec:LEI}, inspired by Gautschi
\cite{rg}, which gives homological constraints on slices of
pseudoholomorphic curves, roughly from the fact that neighborhoods of
these slices have positive energy.

In this paper we use a number of general results about
pseudohomolomorphic curves in $\R\times Y$ which are proved in
\cite{pfh2}, so it may be helpful to have a copy of the latter paper
on hand.

\paragraph{Update} Since the first version of this paper was
distributed, the following related developments have occurred.

In \cite{t3}, we used the techniques of this paper to calculate the
embedded contact homology of $T^3$ in terms of more general
combinatorial chain complexes involving rounding corners of polygons.
Also, \cite{t3} explains how to do such calculations over $\Z$,
instead of over $\Z/2$ as in the present paper.  It follows easily
from \cite{t3} that all homology calculations in the present paper
hold over $\Z$, with all $\Z/2$ summands replaced by $\Z$.

Jabuka and Mark \cite{jm}, using results of Ozsv\'{a}th and Szab\'{o}
\cite{os3}, computed the Ozsv\'{a}th-Szab\'{o} Floer homology $HF^+$
of the mapping tori of some Dehn twists and compositions thereof on
closed surfaces $\Sigma$.  If $d < g(\Sigma)-1$, then the $PFH$ agrees
with $HF^+$ of the mapping torus in those cases where both have been
calculated.  This provides some nontrivial evidence for the
conjectured relation between $PFH$ and Seiberg-Witten Floer homology.

\paragraph{Acknowledgments} 
We thank D.\ Canary, Y.\ Eliashberg, R.\ Gautschi, H.\ Hofer, A.\ 
Ivrii, L.\ Mosher and M.\ Skandera for helpful discussions.

The first author was partially supported by NSF grant DMS-0204681 and
the Alfred P. Sloan Foundation.  The second author was partially
supported by NSF grant DMS-0305825 and a VIGRE fellowship.

\section{Review of periodic Floer homology}
\label{sec:pfh}

Let $\Sigma$ be a compact surface, possibly with boundary, and let
$\omega$ be a symplectic form on $\Sigma$. Let $\phi$ be a
symplectomorphism of $(\Sigma,\omega)$.  Let $Y=Y_\phi$ denote the
mapping torus of $\phi$ as in \eqref{eqn:mappingTorus}, and pick a
``sector'' $h\in H_1(Y)$.  Under the ``monotonicity'' and
``regularity'' assumptions below, we now define the periodic Floer
homology $HP_*(\phi,h)$.  This is the homology of a chain complex
which is generated by ``admissible orbit sets'' and graded by the
relative index $I$, and whose differential counts ``flow lines'', all
defined below.

We use the following notation.  Let $V\to Y$ denote the vertical
tangent bundle of $Y\to S^1$.
Let $[\Sigma]\in H_2(Y)$ denote the
homology class of a fiber, and define the {\em degree\/} $d\eqdef
h\cdot[\Sigma]\in\Z$.
Define the ``index ambiguity class''
\[
c(h) \eqdef c_1(V) + 2h^{\times} \in H^2(Y;\Z),
\]
where $h^{\times}$ denotes the image of $h$ under $H_1(Y)\to
H_1(Y,\partial Y) \simeq H^2(Y;\Z)$.
Let $R$ denote the flow
on $Y$ that increases the $[0,1]$ coordinate in
\eqref{eqn:mappingTorus}.

\paragraph{Orbit sets}
An {\em orbit set\/} is a finite set of pairs $\{(\alpha_i,m_i)\}$
where the $\alpha_i$'s are distinct (nondegenerate) irreducible
periodic orbits of $\phi$, regarded as embedded oriented circles in
$Y$ tangent to $R$, and the $m_i$'s are positive integers, which can
be thought of as ``multiplicities''.  The orbit set
$\{(\alpha_i,m_i)\}$ is {\em admissible\/} if $m_i=1$ whenever
$\alpha_i$ is {\em hyperbolic\/}, i.e.\ the linearized return map of
$\phi$ along $\alpha_i$ has real eigenvalues.  We define the homology
class $[\alpha]\eqdef \sum_im_i[\alpha_i]\in H_1(Y)$, and we let
$\mc{A}(h)$ denote the set of admissible orbit sets $\alpha$ with
$[\alpha]=h$.  We often denote an orbit set by a commutative product
$\alpha_1^{m_1}\cdots\alpha_k^{m_k}$, although the index and
differential defined below are not well-behaved with respect to this
sort of multiplication.

\paragraph{The relative index}
Let $\alpha=\{(\alpha_i,m_i)\}$ and $\beta=\{(\beta_j,n_j)\}$ be orbit
sets with $[\alpha]=[\beta]$.  Let $H_2(Y;\alpha,\beta)$ denote the
set of relative homology classes of 2-chains $W$ in $Y$ with
\[
\partial
W = \sum_im_i\alpha_i - \sum_jn_j\beta_j;
\]
this is an affine space
over $H_2(Y)$.  If $Z\in H_2(Y;\alpha,\beta)$, define the {\em
relative index\/}
\begin{equation}
\label{eqn:RI}
I(\alpha,\beta;Z) \eqdef c_\tau(Z)
+ Q_\tau(Z) + \mu_\tau(\alpha,\beta) \in \Z.
\end{equation}
Here $\tau$ is a homotopy class of symplectic trivialization of $V$
over the $\alpha_i$'s and $\beta_j$'s; $c_\tau(Z)$ denotes the
relative first Chern class of the bundle $V$ over the relative
homology class $Z$ with respect to the boundary trivialization $\tau$;
$Q_\tau(Z)=Q_\tau(Z,Z)$ denotes the relative intersection pairing in
$[0,1]\times Y$; and
\[
\mu_\tau(\alpha,\beta)
\eqdef \sum_i\sum_{k=1}^{m_i}\mu_\tau(\alpha_i^k)
- \sum_j\sum_{k=1}^{n_j}\mu_\tau(\beta_j^k),
\]
where $\mu_\tau(\gamma^k)$ is the Conley-Zehnder index of the $k^{th}$
iterate of $\gamma$.  These notions are explained in detail in
\cite[\S2]{pfh2}.

The relative index has the following basic properties \cite[Prop.\
1.6]{pfh2}.  First, the definition does not depend on $\tau$.  Second,
if $\gamma$ is another homologous orbit set and $W\in
H_2(Y;\beta,\gamma)$, then we have the additivity property
\[
I(\alpha,\beta;Z)+I(\beta,\gamma;W)=I(\alpha,\gamma;Z+W).
\]
Third, if
$\alpha$ and $\beta$ are admissible, then
\begin{equation}
\label{eqn:indexParity}
I(\alpha,\beta;Z)\equiv
\sum_i\epsilon(\alpha_i)-\sum_j\epsilon(\beta_j)\mod 2,
\end{equation}
where $(-1)^{\epsilon}$ denotes the Lefschetz sign, which is $-1$ for
hyperbolic orbits with positive eigenvalues and $+1$ otherwise.
Fourth, we have the ``index ambiguity formula''
\begin{equation}
\label{eqn:indexAmbiguity}
I(\alpha,\beta;Z)-I(\alpha,\beta;W) = \langle Z-W, c(h)\rangle.
\end{equation}

\paragraph{Monotonicity}
The symplectic form $\omega$ on $\Sigma$ induces a symplectic
structure on $V$ which canonically extends to a closed 2-form on $Y$,
which we still denote by $\omega$, with $R\feed\omega=0$.  We say that
$(\phi,h)$ is {\em monotone\/} if
\begin{equation}
\label{eqn:monotone}
[\omega] = \lambda c(h) \in H^2(Y;\R),
\end{equation}
where $\lambda\in\R$.  An elementary calculation in \cite{pfh1} shows
that for any given $h$, if $\Sigma$ is connected and
$\partial\Sigma\neq\emptyset$ or $d\neq g(\Sigma)-1$, then one can
achieve monotonicity by a symplectic isotopy of $\phi$ fixing
$\phi|_{\partial\Sigma}$.

\paragraph{Almost complex structure}
An almost complex structure $J$ on $\R\times Y$ is {\em admissible\/} if:
\begin{itemize}
\item
$J(\partial_s)=R$, where $s$ denotes the $\R$ coordinate.
\item
$J$ is invariant under the obvious $\R$-action on $\R\times Y$.
\item
$J$ is tamed by $\Omega\eqdef\omega+ds\wedge dt$, that is
$\Omega(v,Jv)>0$ for $v\neq 0$.
\end{itemize}
We say that $(\phi,J)$ is {\em $d$-regular\/} if $J$ is admissible and:
\begin{itemize}
\item Each periodic orbit of period $p\le d$ is nondegenerate.
\item
(Local linearity) For each periodic orbit of period $p\le d$, there
are local coordinates on $\Sigma$ near $\gamma$ for which $\phi^p$ is
linear; and there is a tubular neighborhood $N$ of $\gamma$ on which
$J$ sends $V$ to $V$, by a constant matrix on each fiber of the
projection $N\to\gamma$ induced by the projection $Y\to S^1$.  If
$\gamma$ is elliptic, then $J$ is invariant under the flow $R$ in $N$.
\item
Near each component of $\partial\Sigma$ there are local coordinates
$x\in(-\epsilon,0], y\in\R/\Z$ with $J(\partial_x)=\partial_y$ and
$\phi(x,y)=(x,y+\theta)$, where $q\theta\notin\Z$ for all integers
$1\le q\le d$.
\item $J$ is sufficiently close to sending $V$ to $V$.
\end{itemize}
The local linearity condition might not be necessary but is used in
\cite{pfh2,pfh1} to simplify the analysis. The next condition ensures
that flow lines (defined below) do not approach the boundary, by the
maximum principle.  The last condition is used in \cite{pfh2} to
rule out bubbling of closed pseudoholomorphic curves in moduli spaces
of flow lines.

\paragraph{Flow lines}
We now consider $J$-holomorphic curves in $\R\times Y$, where $J$ is
admissible.  The simplest example of a $J$-holomorphic curve in
$\R\times Y$ is $\R\times\gamma$, where $\gamma$ is a periodic orbit;
we call this a {\em trivial cylinder\/}.  More generally, if $C$ is a
$J$-holomorphic curve, an {\em outgoing end\/} at $\gamma$ of
multiplicity $m$ is an end of $C$ asymptotic to $\R^+\times\gamma^m$
as $s\to+\infty$, where $\gamma^m$ denotes a connected $m$-fold cover
of $\gamma$.  {\em Incoming ends\/} are defined analogously with
$s\to-\infty$.

Let $\alpha=\{(\alpha_i,m_i)\}$ and $\beta=\{(\beta_j,n_j)\}$ be orbit
sets with $[\alpha]=[\beta]$.  A {\em flow line\/} from $\alpha$ to
$\beta$ is a $J$-holomorphic curve $C\subset\R\times Y$ such that:
\begin{itemize}
\item $C$ is embedded, except that trivial cylinders may be repeated,
  although these are not allowed to intersect the rest of $C$.
\item
$C$ is a punctured compact Riemann surface and has outgoing
ends at $\alpha_i$ of multiplicity $q_{i,k}$ with
$\sum_kq_{i,k}=m_i$, incoming ends at $\beta_j$ of
multiplicity $q'_{j,k}$ with $\sum_kq'_{j,k}=n_j$, and no other
ends.
\end{itemize}

Let $\M(\alpha,\beta;Z)$ denote the moduli space of flow lines from
$\alpha$ to $\beta$ in the relative homology class $Z$, and
$\M(\alpha,\beta)\eqdef\union_{Z\in H_2(Y;\alpha,\beta)}
\M(\alpha,\beta;Z)$.  Note that $\R$ acts on these moduli spaces by
translation in $\R\times Y$.  If $C\in\M(\alpha,\beta)$, let $\M_C$
denote the component of the moduli space containing $C$.

If $C\in\M(\alpha,\beta)$, we define $I(C)\eqdef I([C]) \eqdef
I(\alpha,\beta;[C])$.  It is shown in \cite[Thm.\ 1.7]{pfh2}, see also
Proposition~\ref{prop:ultimate} below, that if $J$ is generic and
$(\phi,J)$ is $d$-regular, where $d$ is the degree of $\alpha$ and
$\beta$, then $\M_C$ is a manifold and
\begin{equation}
\label{eqn:indexInequality}
\dim(\M_C)\le I(C).
\end{equation}
Moreover, equality holds only if $C$ is ``admissible'', see
\cite[\S4]{pfh2}.  If $C$ has no trivial cylinders, then
``admissible'' means that $\{q_{i,1},q_{i,2},\ldots\} =
\pout(\alpha_i,m_i)$ and
$\{q'_{j,1},q'_{j,2},\ldots\}=\pin(\beta_j,n_j)$, where to any
periodic orbit $\gamma$ and positive integer $m$ there are {\em a priori\/}
associated an ``outgoing partition'' $\pout(\gamma,m)$ and an
``incoming partition'' $\pin(\gamma,m)$ of $m$.  A simple example
which we will need later is that if $\gamma$ is {\em elliptic\/}, i.e.\ the
linearized return map has eigenvalues $e^{\pm 2\pi i\theta}$, and if the
linearized return map is a small (with respect to $d$) clockwise
rotation, then for $m\le d$,
\begin{equation}
\label{eqn:ellipticPartitions}
\pout(\gamma,m)=\{m\},\quad\quad\pin(\gamma,m)=\{1,\ldots,1\}.
\end{equation}
So if $C$ is admissible and has no trivial cylinders, then $C$ can
have at most one outgoing end at $\gamma$, while every incoming end of
$C$ at $\gamma$ has multiplicity one.

\paragraph{The chain complex}
Assume that $(\phi,h)$ is monotone, and that $(\phi,J)$ is $d$-regular
and $J$ is generic.  We now define a chain complex
$(CP_*(\phi,h),\delta)$, whose differential may depend on $J$.  It is
possible to define this over $\Z$, but for simplicity we will work
over $\Z/2$ in the present paper.  The generators of the chain complex
are admissible orbit sets:
\[
CP_*(\phi,h) \eqdef (\Z/2)\{\mc{A}(h)\}.
\]
The differential is defined as follows: if $\alpha\in\mc{A}(h)$, then
\[
\delta\alpha\eqdef
\sum_{\beta\in\mc{A}(h)}
\left(
\sum_{I(\alpha,\beta;Z)=1} \#\frac{\M(\alpha,\beta;Z)}{\R}
\right)\beta,
\]
where `$\#$' denotes the mod 2 count.  By equations
\eqref{eqn:indexAmbiguity} and \eqref{eqn:monotone}, all flow lines on
the right hand side have the same integral of $\omega$, so by the
compactness theorem of \cite[Thm.\ 1.8]{pfh2}, this is a finite sum.
It will be shown in \cite{pfh1} that $\delta^2=0$.  The homology of
this chain complex is the periodic Floer homology $HP_*(\phi,h)$.

By equation \eqref{eqn:indexParity}, $HP_*(\phi,h)$ has a canonical
$\Z/2$ grading, which by \eqref{eqn:indexAmbiguity} lifts
noncanonically to a $\Z/N$ grading, where $N$ is the divisibility
of $c(h)$ in $H^2(Y;\Z)$.

We expect that $HP_*(\phi,h)$ is independent of $J$ and invariant
under symplectic isotopy of $\phi$ fixing $\phi|_{\partial\Sigma}$ and
preserving monotonicity.  Moreover if $\Sigma$ is closed and connected
and $d<g(\Sigma)-1$, then it is conjectured in \cite{pfh1} that
$HP_*(\phi,h)$ agrees with the Seiberg-Witten Floer homology of $Y$
for a spin-c structure determined by $h$.

\paragraph{Duality}
In this paper we mainly discuss positive Dehn twists, but one can
easily deduce corresponding results for negative Dehn twists as
follows.  For any symplectomorphism $\phi$ of $\Sigma$, the
self-diffeomorphism of $\R\times[0,1]\times\Sigma$ sending
$(s,t,x)\mapsto(-s,1-t,x)$ induces a symplectomorphism
$\imath:\R\times Y_{\phi}\to \R\times Y_{\phi^{-1}}$, which sends a
$d$-regular almost complex structure for $\phi$ to one for
$\phi^{-1}$.  Since incoming and outgoing ends of flow lines are
switched, we obtain $HP_*(\phi^{-1},\imath_*h)\simeq
HP^{-*}(\phi,h)$, where the right hand side is the ``periodic Floer
cohomology'' defined using the dual differential.

\paragraph{Generalized flow lines}
In some technical arguments we need to consider {\em generalized flow
lines\/} (GFL's), which are defined like flow lines except that we
drop the embeddedness condition, see \cite[\S9.3]{pfh2}.  We collect
here the definitions and facts about GFL's that we will need.

If $C$ is a GFL, then in the notation from the paragraph on flow
lines, we define
\begin{align}
\nonumber
\mu_\tau^0(C) & \eqdef \sum_i\sum_k\mu_\tau(\alpha_i^{q_{i,k}})
- \sum_j\sum_k\mu_\tau(\beta_j^{q'_{j,k}}),\\
\label{eqn:Ivir}
I^{\op{vir}}(C) &\eqdef -\chi(C)+2c_\tau([C])+\mu_\tau^0(C).
\end{align}
If $C$ is a GFL from $\alpha$ to $\beta$, let $\op{deg}(C)$ denote the
degree of $\alpha$ and $\beta$.

We say that a GFL is a {\em quasi-embedding\/} if it is embedded
except possibly at finitely many points.  A connected GFL is either a
quasi-embedding or a multiple cover.  If $J$ is generic, then it
follows from the index formula of \cite{sch} and a transversality
argument in \cite{pfh2} that if $C$ is a quasi-embedding, then
$I^{\op{vir}}(C)$ equals the dimension of the component of the moduli
space of quasi-embedded GFL's containing $C$.  If $C$ is a flow line,
and if $J$ is generic and $(\phi,J)$ is $\op{deg}(C)$-regular, then
\begin{equation}
\label{eqn:DMV}
\dim(\M_C)\le I^{\op{vir}}(C).
\end{equation}
This is an inequality, because a moduli space of embedded flow lines
can appear as a stratum in a moduli space of GFL's, see
\cite[\S5]{pfh2}.

If $C$ is a quasi-embedded GFL, then by \cite[\S3]{pfh2}, we have the
relative adjunction formula
\begin{equation}
\label{eqn:adj}
c_\tau([C])=\chi(C)+Q_\tau([C])+w_\tau(C)-2\delta(C).
\end{equation}
Here the integer $\delta(C)$ is a weighted count of the singularities
of $C$, which satisfies $\delta(C)\ge 0$ with equality if and only if
$C$ is embedded.  Also, $w_\tau(C)$ is the signed sum of the
writhes of the asymptotic braids of $C$, see \cite[\S3.1]{pfh2}.

It is shown in \cite[\S6]{pfh2} that if $C'$ is a quasi-embedded GFL
without trivial cylinders, and if $(\phi,J)$
is $\op{deg}(C')$-regular, then
\begin{equation}
\label{eqn:WI}
w_\tau(C')\le \mu_\tau(\alpha,\beta)-\mu_\tau^0(C').
\end{equation}
Equations \eqref{eqn:RI}, \eqref{eqn:Ivir}, \eqref{eqn:adj}, and
 \eqref{eqn:WI} imply that
\begin{equation}
\label{eqn:IVI}
I^{\op{vir}}(C')  \le I(C') - 2\delta(C').
\end{equation}
It is also shown in \cite[\S7]{pfh2} that if $C'$ is a quasi-embedded
GFL without trivial cylinders and if $T$ is a union of trivial
cylinders, then
\begin{equation}
\label{eqn:TC}
I(C') \le I(C'\cup T) - 2\#(C'\cap T).
\end{equation}
Here `$\#$' denotes the algebraic intersection number.  Moreover, by
intersection positivity \cite{mcd},
\[
\#(C'\cap T) \ge 0
\]
with equality if and only if $C'\cap T = \emptyset$.

Observe that the inequalities \eqref{eqn:DMV}, \eqref{eqn:IVI}, and
\eqref{eqn:TC} imply the index inequality \eqref{eqn:indexInequality}.
The latter has the following generalization for multiply covered
GFL's.  This is proved in \cite{pfh2}, but the proof below is much
simpler.

\begin{proposition}
\label{prop:ultimate}
Let $C_1,\ldots,C_k$ be distinct, connected, quasi-embedded\break GFL's, and
let $d_1,\ldots,d_k$ be positive integers.  Suppose that $(\phi,J)$ is
$d$-regular, where $d=\sum_id_i\op{deg}(C_i)$.  Then
\begin{equation}
\label{eqn:ultimate}
\sum_{i=1}^k d_i I^{\op{vir}}(C_i) \le
I\left(\sum_{i=1}^kd_i[C_i]\right) - 2\Delta,
\end{equation}
where
\[
\Delta=\sum_{i=1}^k d_i^2\delta(C_i) + \sum_{i<j} d_id_j\#(C_i\cap
C_j).
\]
\end{proposition}

Note that by intersection positivity, $\Delta\ge 0$, with equality if
and only if the $C_i$'s are embedded and disjoint.

\begin{proof}
Suppose first that none of the $C_i$'s is a trivial cylinder, so that
$\R$ acts freely on the $C_i$'s.  Let $\widehat{C}_i$ be a union of
$d_i$ distinct translates of $C_i$, and let
$C'=\bigcup_i\widehat{C}_i$.  Then $C'$ is a quasi-embedded GFL.  Now
the three terms in the inequality \eqref{eqn:IVI} are precisely the
three terms in \eqref{eqn:ultimate}.

If some of the $C_i$'s are trivial cylinders, then use the inequality
\eqref{eqn:TC}, where $C'$ is defined as above for the nontrivial components.
\end{proof}

\begin{corollary}
\label{cor:restrict}
Let $C$ be a GFL with $I(C)=1$.  Assume that $J$ is generic and $(\phi,J)$ is
$\op{deg}(C)$-regular.  Then $C$ has one nontrivial embedded component
$C'$; all other components are covers of trivial cylinders which do
not intersect $C'$.
\end{corollary}

That is, $C$ is a flow line, except that trivial cylinder components
may be nontrivially covered instead of just repeated.

\begin{proof}
We know that $C$ is a union of $d_i$-fold covers of distinct,
  connected, quasi-embedded GFL's $C_i$.  Since $J$ is generic,
  $I^{\op{vir}}(C_i)\ge 0$, with equality if and only if $C_i$ is
  trivial.  So by the inequality \eqref{eqn:ultimate}, there is only
  one nontrivial $C_i$, and this has $d_i=1$ and does not intersect
  the trivial components.
\end{proof}

A {\em {\rm(}$k$ times{\rm)} broken GFL\/} is a sequence $(C_0,\ldots,C_k)$ of
nontrivial GFL's such that for each $i=1,\ldots,k$, the outgoing ends
of $C_i$ are identified with the incoming ends of $C_{i-1}$, such that
two ends that are identified are at the same periodic orbit and with
the same multiplicity.  The broken GFL $(C_0,\ldots,C_k)$ is {\em
  connected\/} if the graph with one vertex for each component of each
$C_i$, and an edge between two vertices when the corresponding
components have ends identified, is connected. A {\em component\/} of
the broken GFL $(C_0,\ldots,C_k)$ is a maximal connected broken GFL
$(C_0',\ldots,C_k')$ such that each $C_i'$ is a union of components of
$C_i$.  By Gromov compactness as in \cite[Lem. 9.8]{pfh2}, any
sequence of nontrivial GFL's of the same topological type with bounded
integral of $\omega$ has a subsequence which converges in an
appropriate sense to a broken GFL.

\section{PFH of a twist on a cylinder}
\label{sec:cylinder}

Fix an integer $P$ and a positive integer $Q$.  Fix real numbers
$X_1\le X_2$, and assume that neither is a rational number with
denominator $\le Q$. Let $\phi_0$ be the cylinder twist from equation
\eqref{eqn:cylinderTwist}.  We identify the mapping torus
\begin{equation}
\label{eqn:cylinderCoordinates}
\begin{array}{rcl}
Y_{\phi_0} & \simeq & S^1\times [X_1,X_2]\times S^1,\\
{[t,(x,y)]} & \mapsto & (t,x,y-xt).
\end{array}
\end{equation}
We now study the periodic Floer homology $HP_*(\phi,h)$, where $\phi$
is a perturbation of $\phi_0$ described below, and the homology class
\begin{equation}
\label{eqn:PQSector}
h
\eqdef
Q[S^1]\times[pt]\times[pt] - P[pt]\times[pt]\times [S^1].
\end{equation}
Note that $(\phi,h)$ is automatically monotone, with $\lambda=0$ in
equation \eqref{eqn:monotone}.  We denote this periodic Floer homology
by $HP_*(X_1,X_2;P,Q)$.

\subsection{Introduction and statement of results}
\label{sec:cylinderResult}

\paragraph{The perturbation}
We always write rational numbers in reduced form $p/q$ with
$\op{gcd}(p,q)=1$ and $q>0$.  For each rational number
$p/q\in[X_1,X_2]$, the map $\phi_0$ has a circle of period $q$
periodic orbits at $x=p/q$.  As is familiar from KAM theory, we can
perform a Hamiltonian perturbation of $\phi_0$ away from the
boundary of $[X_1,X_2]\times S^1$ to obtain a map $\phi$ so that
whenever $q\le Q$, the above circle splits into two periodic orbits
$e_{p/q}$ and $h_{p/q}$.  The orbit $e_{p/q}$ is elliptic,
and the orbit $h_{p/q}$ is hyperbolic with positive eigenvalues.  For
$p/q=0/1$ for example, $\phi_0$ and $\phi$ are the time-one
maps of flows looking like this:
\begin{center}
{\small
\psfraga <-3pt, 1pt> {h}{$h$}
\psfraga <-3pt, 1pt> {e}{$e$}
\psfraga <-3pt, 0pt> {x}{$x$}
\psfraga <-3pt, 0pt> {y}{$y$}
\psfraga <-3pt, 1pt> {0/1}{\tiny$0/1$}
\includegraphics[width=9cm]{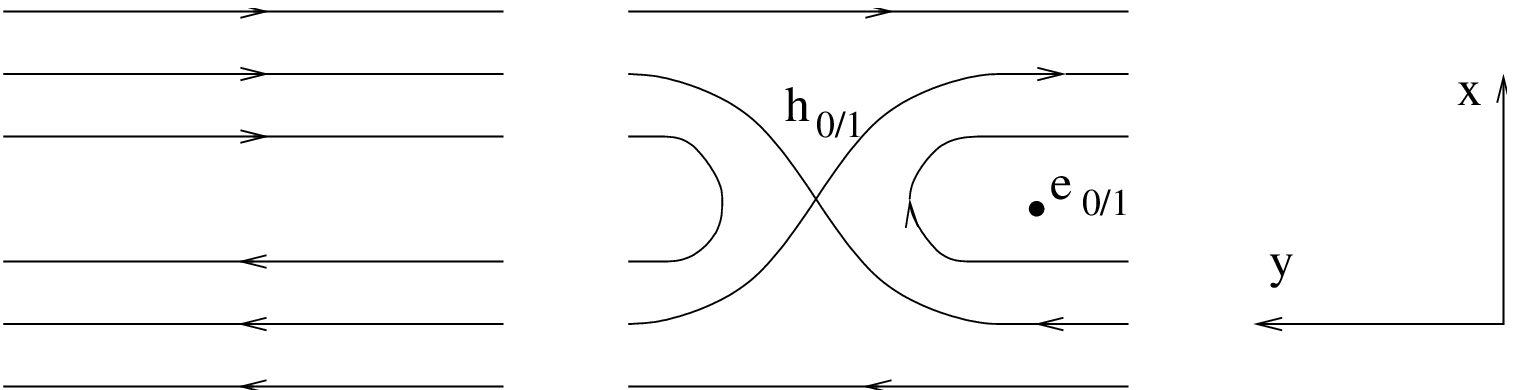}}
\end{center}

Throughout \S\ref{sec:cylinder}, we assume unless otherwise stated that:
\begin{itemize}
\item $\phi$ as above is chosen so that the $e_{p/q}$'s and
  $h_{p/q}$'s are the only irreducible periodic orbits of period $\le
  Q$; these are nondegenerate; and the linearized return map of the
  elliptic orbit $e_{p/q}$ has eigenvalues $e^{\pm 2\pi i\theta}$ with
  $0 < \theta < 1/Q$.
\item
  $\phi$ agrees with $\phi_0$ whenever $x$ is not within distance
  $\varepsilon$ of a rational number of denominator $\le Q$,
  where $\varepsilon$ is sufficiently small with respect to $Q$.
\item $J$ is a generic admissible almost complex structure on
  $\R\times Y_\phi$, such that the pair $(\phi,J)$ satisfies the local
  linearity condition in \S\ref{sec:pfh} with $d=Q$.
\end{itemize}
(In this example the last two conditions in the definition of
$d$-regular are not needed to define PFH: bubbling of closed
pseudoholomorphic curves in $\R\times Y$ cannot happen here because
$\Omega$ is zero on $H_2(\R\times Y_{\phi})$, and flow
lines cannot approach the boundary by Lemma~\ref{lem:slice} below.  In
some arguments we will also drop the local linearity assumption.)

We denote the corresponding chain complex by $CP_*(X_1,X_2;P,Q)$.  We
will see below that the homology of this chain complex, as well as
most of the differential on it, does not depend on the choice of
$\phi$ and $J$ as above.

\paragraph{The generators}
By \eqref{eqn:PQSector}, our chain complex is generated by admissible
products of $e_{p/q}$'s and $h_{p/q}$'s with total numerator $P$ and
total denominator $Q$.  Recall that ``admissible'' means that no
$h_{p/q}$ factor may be repeated.

To each generator $\alpha$ we associate a left-turning convex
polygonal path $\mc{P}(\alpha)$ in the plane as follows.  Write
$\alpha=\gamma_1\cdots\gamma_k$ where 
$\gamma_i=e_{p_i/q_i}$ or $\gamma_i=h_{p_i/q_i}$,
and $\sum_{i=1}^k(p_i,q_i)=(P,Q)$.  We order the factors so that
$p_i/q_i\ge p_j/q_j$ for $i<j$.  For $j=0,\ldots,k$, let
$w_j=\sum_{i=0}^j(p_i,q_i)$, and define $\mc{P}(\alpha)$ to be the
convex path in the plane consisting of straight line segments between
the points $w_{j-1}$ and $w_j$ for $j=1,\ldots,k$, oriented so that
the origin is the initial endpoint.

\paragraph{The homology}
Suppose $P/Q\in[X_1,X_2]$.  We can uniquely write $(P,Q)=v_1+v_2$,
where the vector $v_i$ is in the upper half plane and has slope
$X_i^{-1}$.  Let
$Z=Z(X_1,X_2;P,Q)\subset \R^2$
denote the parallelogram with vertices $0$, $v_1$, $v_2$, and $(P,Q)$.
We define
\[
E(X_1,X_2;P,Q)\in CP_*(X_1,X_2;P,Q)
\]
to be the unique generator $E$ such that $E$ is a product of elliptic
factors, and the path $\mc{P}(E)$ traverses the right half of the
boundary of the convex hull of the set of lattice points in the
parallelogram $Z$.
For example, if $-1/4<X_1<-1/5$ and
$4/3<X_2<3/2$, then
\[
E\left(X_1,X_2;4,11\right)=
e_{4/3}e_1e_0^2e_{-1/5},
\]
as shown by the following picture:

\begin{center}
{\small
\psfrag{p}{$p$}
\psfrag{q}{$q$}
\includegraphics[width=7cm]{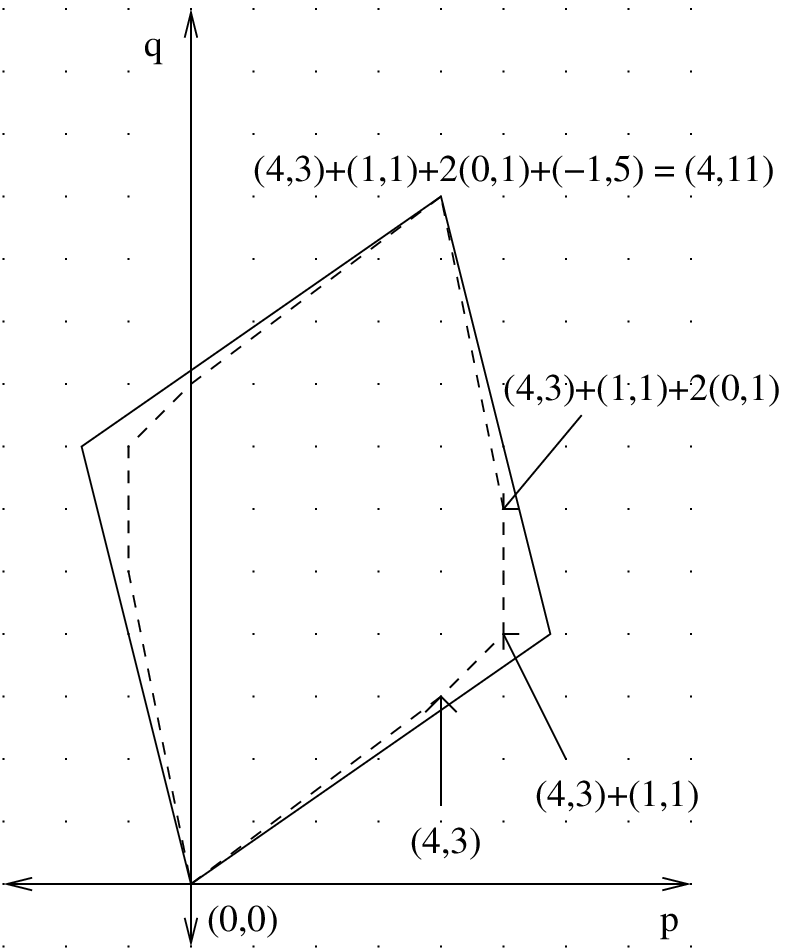}}
\end{center}

\begin{theorem}[PFH of a twist on a cylinder]
\label{thm:cylinder}
\[
HP_*(X_1,X_2;P,Q) \simeq \left\{ \begin{array}{cl}
H_*(S^1;\Z/2) & \mbox{if\qua $P/Q\in[X_1,X_2]$,}\\
0 & \mbox{if\qua $P/Q\notin[X_1,X_2]$.}
\end{array}
\right.
\]
If $P/Q\in [X_1,X_2]$, then $HP_*(X_1,X_2;P,Q)$ is generated by:
\begin{itemize}
\item
$E(X_1,X_2;P,Q)$, and
\item
any generator obtained from $E(X_1,X_2;P,Q)$ by replacing one of the
 $e_{p/q}$ factors by $h_{p/q}$; all such generators are all homologous.
\end{itemize}
\end{theorem}

\paragraph{Combinatorial formulas}  Theorem~\ref{thm:cylinder} is proved in
\S\ref{sec:cylinderSpectral}.  First, in
\S\ref{sec:index}--\S\ref{sec:nonvanishing}, we derive
combinatorial formulas for the relative index and most of the
differential $\delta$ in the chain complex $CP_*(X_1,X_2;P,Q)$.  To
state these, let $\alpha$ and $\beta$ be generators.  We will
see in \S\ref{sec:index} that the relative index
$I(\alpha,\beta;Z)\in\Z$ does not depend on $Z$, so we write it as
$I(\alpha,\beta)$, thus obtaining a
relative $\Z$-grading on $CP_*(X_1,X_2;P,Q)$.  We then have:

\begin{proposition}[Index formula]
\label{prop:area}
If there are $a$ elliptic factors in $\alpha$ and $b$ elliptic factors
in $\beta$, then
\[
I(\alpha,\beta)
=
b-a+2\int_{\mc{P}(\beta)-\mc{P}(\alpha)}p\,dq.
\]
\end{proposition}

\begin{remark}
\label{rem:pick}
Here $\int_{\mc{P}(\beta)-\mc{P}(\alpha)}p\,dq$ denotes the signed
area enclosed by the paths $-\mc{P}(\alpha)$ and $\mc{P}(\beta)$.
Pick's theorem confirms that this is a half-integer.  If
$\mc{P}(\alpha)$ never crosses to the right of $\mc{P}(\beta)$ and if
$\alpha$ and $\beta$ contain no hyperbolic factors, then
$I(\alpha,\beta)$ equals twice the number of lattice points in the
region enclosed by $\mc{P}(\alpha)$ and $\mc{P}(\beta)$, not including
the lattice points in $\mc{P}(\alpha)$.
\end{remark}

If $\mc{P}_1$ and $\mc{P}_2$ are paths in the $(p,q)$-plane, let
$\mc{P}_1\mc{P}_2$ denote the concatenation of $\mc{P}_1$ with the
appropriate translate of $\mc{P}_2$.

\begin{definition}
\label{def:rounding}
We say that $\alpha$ is obtained from $\beta$ by {\em rounding a
corner\/} if there exist orbit sets $\gamma_1,\gamma_2,\alpha',\beta'$
such that:
\begin{itemize}
\item[(a)]
$\alpha=\gamma_1\alpha'\gamma_2$ and $\beta=\gamma_1\beta'\gamma_2$,
where $\mc{P}(\alpha)=\mc{P}(\gamma_1)\mc{P}(\alpha')\mc{P}(\gamma_2)$ and
$\mc{P}(\beta)=\mc{P}(\gamma_1)\mc{P}(\beta')\mc{P}(\gamma_2)$.
\item[(b)]
The path $\mc{P}(\alpha')$ does not cross to the right of
$\mc{P}(\beta')$, and only intersects it at the two endpoints.
\item[(c)]
There
are no lattice points in between $\mc{P}(\alpha')$ and
$\mc{P}(\beta')$.
\item[(d)]
$\beta'$ has two factors; either $\beta'$ has one hyperbolic factor
and $\alpha'$ has none, or $\beta'$ has two hyperbolic factors and
$\alpha'$ has one.
\end{itemize}
We say that $\alpha$ is obtained from $\beta$ by {\em double
rounding\/} if (a), (b), and (c) above hold, together with:
\begin{itemize}
\item[(${\rm d}'$)]
$\beta'$ has three factors, all hyperbolic, and $\alpha'$ has no
hyperbolic factors.
\end{itemize}
\end{definition}

\begin{theorem}[The differential]
Let $\alpha,\beta$ be generators of $CP_*(X_1,X_2;P,Q)$.
Then for any $\phi$ and $J$ as in \S\ref{sec:cylinderResult}:
\label{thm:rounding}
\begin{itemize}
\item[\rm(a)]
If $\langle\delta\alpha,\beta\rangle=1$, then $\alpha$ is obtained
from $\beta$ by rounding a corner or by double rounding.
\item[\rm(b)]
If $\alpha$ is obtained from $\beta$ by rounding a corner, then
$\langle\delta\alpha,\beta\rangle=1$.
\end{itemize}
\end{theorem}
Here $\langle\delta\alpha,\beta\rangle\in\Z/2$ denotes the coefficient
of $\beta$ in $\delta\alpha$.

\begin{remark}
\label{rem:MorseBott}
If we drop the local linearity requirement (in which case more
analytic work would be needed to show that PFH is well-defined) and
consider $(\phi,J)$ close in an appropriate norm to $(\phi_0,J_0)$,
where $J_0$ is a certain $S^1$-invariant almost complex structure on
$Y_{\phi_0}$, then we can arrange that there are no flow lines from
$\alpha$ to $\beta$ when $\alpha$ is obtained from $\beta$ by double
rounding, see Appendix~\ref{app:MorseBott}.  In any case, we will see
in \S\ref{sec:cylinderSpectral} that differential coefficients
involving double rounding have no effect on the homology.
\end{remark}

\subsection{Calculating the relative index}
\label{sec:index}

Let $\alpha$ and $\beta$ be generators of $CP_*(X_1,X_2;P,Q)$ and let
$Z\in H_2(Y;\alpha,\beta)$. We now prove Proposition~\ref{prop:area},
computing the relative index $I(\alpha,\beta;Z)$.

The bundle $V\simeq T([X_1,X_2]\times S^1)$ over $Y_\phi\simeq
S^1\times [X_1,X_2]\times S^1$ has a canonical trivialization up to
homotopy, giving rise to a natural trivialization $\tau$ over each
periodic orbit, which we use for the rest of \S\ref{sec:index}.

Because the trivialization $\tau$ comes from a global trivialization
of $V$,
\begin{equation}
\label{eqn:zeroChern}
c_\tau(Z) = 0.
\end{equation}
We now compute the Conley-Zehnder indices as in \cite[\S2.3]{pfh2}.
The elliptic orbit $e_{p/q}$ has slightly negative monodromy angle
$\theta$ with respect to $\tau$, so
\begin{equation}
\label{eqn:ellInd}
\mu_\tau(e_{p/q}^k) = 2 \lfloor k \cdot \theta \rfloor + 1 = -1
\end{equation}
for $kq\le Q$.  The hyperbolic orbit $h_{p/q}$ has positive
eigenvalues, and the eigenspaces do not rotate with respect to $\tau$,
so
\begin{equation}
\label{eqn:hypInd}
\mu_\tau(h_{p/q}) = 0.
\end{equation}
We next observe that $Q_\tau(Z)$ depends only on $\alpha$ and
$\beta$, by \cite[Lem.\ 2.5]{pfh2}, since $H_2(Y)$ is generated by a
$(y,t)$-torus, which has algebraic intersection number zero with every
periodic orbit.  So denote the integer $Q_\tau(Z)$ by
$Q_\tau(\alpha,\beta)$; by \eqref{eqn:zeroChern}, we can likewise
write $I(\alpha,\beta)\eqdef I(\alpha,\beta;Z)$.

\begin{lemma}
\label{lem:QFormula}
If $\alpha$ and $\beta$ are generators of $CP_*(X_1,X_2;P,Q)$ then
\begin{equation}
\label{eqn:QFormula}
Q_\tau(\alpha,\beta)=2\int_{\mc{P}(\beta)-\mc{P}(\alpha)}p\,dq.
\end{equation}
\end{lemma}

\begin{proof}
Write $\alpha=\gamma_1\cdots\gamma_k$ and
$\beta=\gamma_1'\cdots\gamma_l'$, where $\gamma_i=e_{p_i/q_i}$ or
$\gamma_i=h_{p_i/q_i}$, and $\gamma_i'=e_{p_i'/q_i'}$ or
$\gamma_i'=h_{p_i'/q_i'}$.  We order the factors so that $p_i/q_i\ge
p_j/q_j$ and $p_i'/q_i'\ge p_j'/q_j'$ for $i<j$.  We claim that
\begin{equation}
\label{eqn:detFormula}
Q_\tau(\alpha,\beta)  = 
- \sum_{i<j} \det
\left(
\begin{matrix} 
p_i & p_j \\
q_i & q_j
\end{matrix}
\right) 
+  \sum_{i<j} \det
\left(
\begin{matrix} 
p'_i & p'_j \\
q'_i & q'_j
\end{matrix}
\right).
\end{equation}
By \cite[Eq.\ (13)]{pfh2}, if $S\subset[0,1]\times Y$ is an embedded
(except at the boundary) representative of a relative homology class
$Z\in H_2(Y;\alpha,\beta)$, then
\begin{equation}
\label{eqn:computeQ}
Q_\tau(\alpha,\beta)=c_1(N,\tau)-w_\tau(S),
\end{equation}
where $N$ denotes the normal bundle to $S$ and $c_1(N,\tau)$ its
relative first Chern class.  To prove \eqref{eqn:detFormula}, we apply
\eqref{eqn:computeQ} to the following surface $S$.  We regard $S$ as a
movie of curves in $Y$ parametrized by $s\in[0,1]$.  For $s$ close to
$1$, the slice $S\cap (\{s\}\times Y)$ consists of one circle in the
$x=p/q$ torus for each factor $\gamma_{p/q}$ in $\alpha$, parallel to
the periodic orbits $e_{p/q}$ and $h_{p/q}$.  As $s$ decreases to
$2/3$, we translate all of the circles in the $x$ direction into a
single $(y,t)$-torus.  Around $s=2/3$, we perform ``negative
surgeries'' so that at $s=1/2$, we have an embedded union of circles
in a single $(y,t)$-torus.  Near a negative surgery, the $s>2/3$ and
$s<2/3$ slices look like the left and right sides of the following
picture:
\begin{center}
{\small
\psfraga <0pt,2pt> {x}{$x$}
\psfrag{y}{$y$}
\psfraga <-2pt,0pt> {t}{$t$}
\includegraphics[width=7cm]{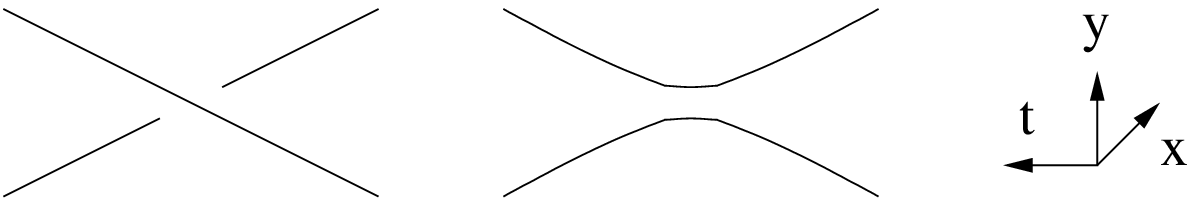}}
\end{center}
For $s$ between $1/2$ and $0$ we perform an opposite process, doing
``positive surgeries'' at $s=1/3$, thereby obtaining, for $s$ close to
zero, one circle in the $x=p/q$ torus for each factor $\gamma_{p/q}$ in
$\beta$.  We can assume that on the complement of the surgery points,
our circles never point in the $y$ direction.

By construction, $w_\tau(S)=0$.  To compute $c_1(N,\tau)$, let $\psi
\eqdef\pi_N\partial_y\in\Gamma(N)$.  Since $\psi$ has winding
number zero with respect to $\tau$ on $\partial S$, the integer
$c_1(N,\tau)$ equals the algebraic count of zeroes of $\psi$.  The section
$\psi$ vanishes only at the surgery points.  The number of negative
surgeries equals the first sum of determinants (without the minus
sign) on the right side of equation
\eqref{eqn:detFormula}, and the number of positive surgeries equals
the second sum.  It is an exercise to check that $\psi$ has $-1$
zero at each negative surgery point and $+1$ zero at each positive
surgery point.  This proves \eqref{eqn:detFormula}.

To deduce \eqref{eqn:QFormula}, observe that the right
hand sides of \eqref{eqn:QFormula} and \eqref{eqn:detFormula} change
by the same amount if we replace two consecutive edges in
$\mc{P}(\alpha)$ or $\mc{P}(\beta)$ by their sum.  This inductively
reduces to the trivial case where $\mc{P}(\alpha)$ and $\mc{P}(\beta)$
are both straight lines.
\end{proof}

\begin{proof}[Proof of Proposition~\ref{prop:area}]  This follows immediately
from equations \eqref{eqn:zeroChern},\break \eqref{eqn:ellInd},
\eqref{eqn:hypInd} and \eqref{eqn:QFormula}.
\end{proof}

\subsection{Some constraints on index 1 flow lines}
\label{sec:constraints}

\begin{lemma}[Complexity]
\label{lem:complexity}
Let $\alpha,\beta\in CP_*(X_1,X_2;P,Q)$ be generators with
$I(\alpha,\beta)=1$, let $C\in\M(\alpha,\beta)$, and
$C'\eqdef C \setminus \{\mbox{trivial cylinders}\}$. Then:
\begin{itemize}
\item[\rm(a)]
$C'$ is connected and has genus zero.
\item[\rm(b)] If $e_-(C')$ denotes the number of incoming elliptic ends
of $C'$, and $h(C')$ denotes the number of hyperbolic ends of $C'$,
then
\[
2e_-(C') + h(C') = 3.
\]
\item[\rm(c)] Any incoming elliptic end of $C'$ has multiplicity one.
  By contrast, $C'$ has only one outgoing end at any given elliptic
  orbit.
\end{itemize}
\end{lemma}

\begin{proof}
  We have $1\le \dim(\M_{C'}) \le I^{\op{vir}}(C')\le I(C') \le
  I(C)=1$, since $\R$ acts nontrivially on $\M_{C'}$ by translation in
  $\R\times Y$, and using the inequalities \eqref{eqn:DMV},
  \eqref{eqn:IVI} and \eqref{eqn:TC}.  Then $C'$ is connected, since
  otherwise $\dim(\M_{C'})\ge 2$.
  
  Let $e_+(C')$ denote the number of outgoing elliptic ends of $C'$.
  By \eqref{eqn:Ivir}, \eqref{eqn:zeroChern}, \eqref{eqn:ellInd}, and
  \eqref{eqn:hypInd},
\[
I^{\op{vir}}(C')=-\chi(C')-e_+(C')+e_-(C').
\]
Putting in $I^{\op{vir}}=1$ and $\chi=2-2g-e_--e_+-h$, we obtain
\begin{equation}
\label{eqn:pre3Dollar}
2g(C')+2e_-(C')+h(C')=3.
\end{equation}
If $g(C')>0$, then \eqref{eqn:pre3Dollar} implies that
$C'\in\M(\alpha',\beta')$ where $\alpha'$ is a product of elliptic
factors and $\beta'=h_{p/q}$ for some $p/q$.  Then by
Proposition~\ref{prop:area}, $I(\alpha',\beta')<0$, because the area
term is $\le 0$ by convexity.  This contradicts $I(C')=1$.

The above proves parts (a) and (b) of the lemma.  Part (c) follows
from the discussion of incoming and outgoing partitions in
\S\ref{sec:pfh}, see equation \eqref{eqn:ellipticPartitions}.
\end{proof}

\begin{lemma}[Trivial cylinders]
\label{lem:trivCyl}
Let $\alpha\gamma,\beta\gamma\in CP_*(X_1,X_2;P,Q)$ be generators,
where $\alpha$ and $\beta$ have no periodic orbits in common.  Assume
$
I(\alpha,\beta)=I(\alpha\gamma,\beta\gamma)=1$.  Then attaching
trivial cylinders over $\gamma$ gives a bijection
\[
\M(\alpha,\beta)\simeq\M(\alpha\gamma,\beta\gamma).
\]
\end{lemma}

\begin{proof}
  By \eqref{eqn:TC} and intersection positivity, if
  $I(C')=I(C'\cup T)$ then $C'\cap T=\emptyset$.  So since
  $I(\alpha,\beta)=I(\alpha\gamma,\beta\gamma)$, it follows that
  attaching trivial cylinders gives a well defined map
  $\M(\alpha,\beta) \to \M(\alpha\gamma,\beta\gamma)$, and
  this is clearly injective.
  
  Now this map must be surjective.  For suppose there exists
  $C\in\M(\alpha\gamma,\beta\gamma)$ which does not contain trivial
  cylinders over all of $\gamma$.  Let $C'$ denote the nontrivial
  component of $C$.  Then $C'$ has both an incoming end and an
  outgoing end at some orbit $\rho$ in $\gamma$.  If $\rho$ is
  hyperbolic, then $\rho$ is not in the generators $\alpha$ and
  $\beta$ (since $\alpha\gamma$ and $\beta\gamma$ are admissible orbit
  sets).  If $\rho$ is elliptic, then by
  Lemma~\ref{lem:complexity}(c), $C'$ has an incoming end of
  multiplicity one there.  Either way, ends of $C'$ at $\gamma$ alone
  contribute at least $2$ to the left hand side of equation
  \eqref{eqn:pre3Dollar}.  Hence Lemma~\ref{lem:complexity}(b) implies
  that $\beta=h_{p/q}$ for some $p/q$, and $\alpha$ is a product of
  elliptic factors.  By Proposition~\ref{prop:area},
  $I(\alpha,\beta)<0$, which is a contradiction.
\end{proof}

\begin{remark}
  The analogue of Lemma~\ref{lem:trivCyl} for a single Dehn twist on a
  torus is false.  We will see in \S\ref{sec:degeneration} that in the
  notation of that section, there exists $C\in\mc{M}(e^2,eh;Z)$ with
  no trivial cylinders even though
\[
I(e^2,eh;Z)=I(e,h;Z)=1.
\]
\end{remark}

\subsection{The local energy inequality}
\label{sec:LEI}

\paragraph{The general case}
Let $\phi$ be a general symplectomorphism of $\Sigma$.  Suppose
$\xi\subset\Sigma$ is an invariant circle, i.e.\ $\phi(\xi)=\xi$, such
that $\phi|_\xi$ is smoothly conjugate to a rotation.  Let $T\eqdef
Y_{\phi|_\xi}$ and let $N\subset\Sigma$ be a tubular neighborhood of
$\xi$.  We choose local coordinates $x\in(-\epsilon,\epsilon)$ and
$y,t\in\R/\Z$ on $Y_{\phi|_N}$ such that $T=\{x=0\}$, the mapping
torus flow at $x=0$ is $R=\partial_t-\theta\partial_y$ with
$\theta\in\R$ fixed, and $\omega(\partial_x,\partial_y)>0$.  Let $J$
be an admissible almost complex structure on $\R\times Y_\phi$.

Let $C$ be a generalized flow line transverse to $\R\times T$.  Orient
$C\cap(\R\times T)$ so that if $p\in C\cap(\R\times T)$ and $\{v,w\}$
is an oriented basis of $T_pC$ with $v\in T_p(C\cap(\R\times T))$
positively oriented, then $w$ has negative $\partial_x$ component.

\begin{lemma}
\label{lem:slice}
The homology class of the slice
\[
(p,q)\eqdef [C\cap(\R\times T)] \in H_1(S^1_y\times S^1_t)
\]
satisfies $p+\theta q\ge 0$, with equality iff $C\cap(\R\times
T)=\emptyset$.
\end{lemma}

\begin{proof}
Each component of $C\cap(\R\times T)$ can be described by an oriented
parametrized curve $\gamma(\tau)$.  Let $a(y,t)$ denote the component
of $\partial_x$ in $J\partial_y$ at $x=0$.  Since $J$ is
$\Omega$-tame, $a(y,t)<0$.  At $x=0$, since
$J(\partial_t-\theta\partial_y)=-\partial_s$, the component of
$\partial_x$ in $J\partial_s$ is zero, and the component of
$\partial_x$ in $J\partial_t$ is $\theta a(y,t)$.  Since $C$ is
$J$-holomorphic, $J(\gamma'(\tau))\in TC$, and by transversality and
our sign conventions, this has negative $\partial_x$ component.  That
is, if $\gamma(\tau)=(s(\tau),y(\tau),t(\tau))$, then
\[
\frac{dy}{d\tau} + \theta\frac{dt}{d\tau}>0.
\]
This inequality can also be understood as intersection positivity of
$C$ with the $J$-holomorphic foliation of $\R\times T$ by $\R$ cross
the mapping torus flow.  Integrating this inequality over $\tau$, we
conclude that the homology class $(p',q')$ of $\gamma$ satisfies
$p'+\theta q'>0$.
\end{proof}

\paragraph{The cylinder twist case} For a Dehn twist on a cylinder,
the local energy inequality of Lemma~\ref{lem:slice} has the following
interpretation.  

\begin{proposition}
\label{prop:leftRight}
Let $\alpha$ and $\beta$ be generators of $CP_*(X_1,X_2;P,Q)$.
Suppose that the path $\mc{P}(\alpha)$ crosses to the right of
$\mc{P}(\beta)$.  Then for any $\phi$ as in \S\ref{sec:cylinderResult}
and any admissible $J$, there are no flow lines from $\alpha$ to
$\beta$.
\end{proposition}

\begin{proof}
  Suppose there exist flow lines from $\alpha$ to $\beta$ for
  arbitrarily small values of the constant $\varepsilon$ in
  \S\ref{sec:cylinderResult}.  Let $p_i,q_i,p_j',q_j'$ be defined as
  in the proof of Lemma~\ref{lem:QFormula}.  We first claim that for
  each $x_0\in\R$, if $\varepsilon$ is sufficiently small then
\begin{equation}
\label{eqn:LEI}
0 \le
\sum_{p_i/q_i > x_0} ( q_i x_0 - p_i) - 
\sum_{p'_j/q'_j > x_0}  (q'_j x_0 - p'_j).
\end{equation}
To prove \eqref{eqn:LEI}, by continuity we may assume that $x_0$ is
not a rational number of denominator $\le Q$.  We can then assume that
$\varepsilon$ is sufficiently small that $\phi$ agrees with $\phi_0$ for
$x$ near $x_0$.  By continuity again, we may assume that $C$
intersects the locus $\{x=x_0\}$ transversely.
Then \eqref{eqn:LEI} follows from Lemma~\ref{lem:slice} since
\begin{equation}
\label{eqn:sliceClass}
[C\cap\{x=x_0\}]
=
\sum_{p_i/q_i > x_0} ( -p_i,q_i) - 
\sum_{p'_j/q'_j > x_0}  (-p'_j, q'_j) \in H_1(S^1_y \times S^1_t).
\end{equation}
Now suppose $(a,b)$ and $(c,d)$ are two intersections of the paths
$\mc{P}(\alpha)$ and $\mc{P}(\beta)$, between which $\mc{P}(\alpha)$
is to the right of $\mc{P}(\beta)$.  Since both paths move in the
positive $q$ direction we may assume that $b<d$.  Let $(p,q)$ and
$(p',q')$ be points on $\mc{P}(\alpha)$ and $\mc{P}(\beta)$
respectively between $(a,b)$ and $(c,d)$ which are as far as possible
from the line through $(a,b)$ and $(c,d)$.  Then $(p,q)$ is farther
from this line than $(p',q')$, so
\begin{equation}
\label{eqn:furtherToRight}
\det\begin{pmatrix}p-p' & c-a\\ q-q' & d-b\end{pmatrix} > 0.
\end{equation}
If we let
$x_0 \eqdef (c-a)/(d-b)$
then by equation \eqref{eqn:LEI} we have
\[
0\le (qx_0-p) - (q'x_0-p').
\]
This inequality contradicts $b<d$ and \eqref{eqn:furtherToRight}.
\end{proof}

\begin{lemma}
\label{lem:inclusion}
If $[X_1,X_2]\subset [X_1',X_2']$, then for compatible choices of $J$
and $\phi$ as in \S\ref{sec:cylinderResult}, the differential $\delta$
commutes with the inclusion
\[
CP_*(X_1,X_2;P,Q)\longrightarrow
CP_*(X_1',X_2';P,Q).
\]
\end{lemma}

\begin{proof}
If $\alpha,\beta\in CP_*(X_1,X_2;P,Q)$, then any flow line from
$\alpha$ to $\beta$ in $S^1\times [X_1',X_2']\times S^1$ is supported
in $S^1\times[X_1,X_2]\times S^1$, by Lemma~\ref{lem:slice}.
\end{proof}

\subsection{Vanishing of some differential coefficients}
\label{sec:vanishing}

\begin{lemma}
\label{lem:MorseBott}
If $\lambda q \le Q$, then for any $J$ and $\phi$ as in
\S\ref{sec:cylinderResult},
\[
\left\langle\delta\left(
e_{p/q}^{\lambda -1}h_{p/q}
\right),
e_{p/q}^\lambda
\right\rangle=0.
\]
\end{lemma}

\begin{proof}
The trivial cylinder lemma \ref{lem:trivCyl} shows that
\[
\left\langle\delta\left(e_{p/q}^{\lambda-1}h_{p/q}\right),e_{p/q}^{\lambda}
\right\rangle
=
\left\langle\delta h_{p/q},e_{p/q}\right\rangle.
\]
By Lemma~\ref{lem:complexity}(a), any flow line from $h_{p/q}$ to
$e_{p/q}$ is a cylinder.  Now for the unperturbed map $\phi_0$ with an
almost complex structure such as $J_0$ in
Appendix~\ref{app:MorseBott}, any pseudoholomorphic cylinder with
incoming and outgoing ends at $x=p/q$ is trivial, as in \cite[Prop.\
9.1]{pfh2}.  Then a standard argument in Morse-Bott theory shows that
for $(\phi',J')$ close to $(\phi_0,J_0)$, there are two cylinders from
$h_{p/q}$ to $e_{p/q}$.  These cylinders are both embedded by the
adjunction formula \eqref{eqn:adj}, because all terms in
\eqref{eqn:adj} other than $2\delta(C)$ are automatically zero.
As we deform $(\phi',J')$ to $(\phi,J)$, no broken GFL's from
$h_{p/q}$ to $e_{p/q}$ appear during the deformation by
Lemma~\ref{lem:slice}.  It follows by Gromov compactness as in
\cite[Lem.\ 9.8]{pfh2} that the mod 2 count of such cylinders remains zero.
\end{proof}

\begin{proof}[Proof of Theorem~\ref{thm:rounding}(a)]
Suppose $\langle\delta\alpha,\beta\rangle\neq 0$.

Without loss of generality, $\mc{P}(\alpha)$ and
$\mc{P}(\beta)$ have no initial edges in common.  To see this, suppose
that $\alpha$ and $\beta$ both contain factors of $e_{p/q}$ or
$h_{p/q}$ but do not contain factors $e_{p'/q'}$ or $h_{p'/q'}$ with
$p'/q'>p/q$.  We claim that $\alpha$ and $\beta$ have a common factor
of $e_{p/q}$ or $h_{p/q}$, so that we can remove it from both by the
trivial cylinder lemma \ref{lem:trivCyl}.  Otherwise, either: (1)
$\alpha$ contains an $h_{p/q}$ factor and $\beta$ contains an
$e_{p/q}^{\lambda}$ factor, or (2) $\alpha$ contains an $e_{p/q}^\lambda$
factor and $\beta$ contains an $h_{p/q}$ factor.  In case (1), we must
have $\lambda=1$ by Lemma~\ref{lem:complexity}(b).  By
Lemma~\ref{lem:slice}, if $C\in\mc{M}(\alpha,\beta)$, then
$C\cap\{x=p/q-\varepsilon\}=\emptyset$, so
\[
\mc{M}(\alpha,\beta)=
\mc{M}(\alpha/h_{p/q},\beta/e_{p/q})\times\mc{M}(h_{p/q},e_{p/q}),
\]
so $\langle\delta\alpha,\beta\rangle=0$ by Lemma~\ref{lem:MorseBott}.
In case (2), we have $\lambda=1$ by Proposition~\ref{prop:leftRight}.
Then as in case (1),
\[
\mc{M}(\alpha,\beta)=
\mc{M}(\alpha/e_{p/q},\beta/h_{p/q})\times\mc{M}(e_{p/q},h_{p/q}).
\]
But $\mc{M}(e_{p/q},h_{p/q})=\emptyset$, since
$I(e_{p/q},h_{p/q})=-1$, so $\langle\delta\alpha,\beta\rangle=0$.

Likewise, WLOG the paths $\mc{P}(\alpha)$ and $\mc{P}(\beta)$ have no
final edges in common.

By Proposition~\ref{prop:leftRight}, $\mc{P}(\alpha)$ does not cross to
the right of $\mc{P}(\beta)$.  Moreover, the paths $\mc{P}(\alpha)$
and $\mc{P}(\beta)$ do not intersect except at their endpoints.  For
if the paths intersect elsewhere, then the intersections must be
lattice points or segments bounded by lattice points.  These intersections cut
$\mc{P}(\beta)$ into at least two pieces.  By convexity, each piece
has at least two edges, and the initial and final edges of each piece
do not have parallel edges in $\mc{P}(\alpha)$.  It follows that,
aside from possible trivial cylinders, flow lines from $\alpha$ to
$\beta$ have at least four incoming ends, contradicting
Lemma~\ref{lem:complexity}(b).

Thus conditions (a) and (b) in Definition~\ref{def:rounding} hold.
Since $I(\alpha,\beta)=1$, conditions (c) and either (d) or (d$'$) in
Definition~\ref{def:rounding} follow from the index formula of
Proposition~\ref{prop:area} and Remark~\ref{rem:pick}.
\end{proof}

\subsection{Invariance of rounding coefficients}
\label{sec:IRC}

Suppose that $\alpha$ is obtained from $\beta$ by rounding a corner,
and that $\alpha$ and $\beta$ have no periodic orbits in common.  By
Lemma~\ref{lem:complexity}, the flow lines $C\in\mc{M}(\alpha,\beta)$
counted by the differential coefficient
$\langle\delta\alpha,\beta\rangle$ satisfy the following conditions:
\begin{itemize}
\item[\rm(i)] $C$ is connected and has genus zero.
\item[\rm(ii)]
For each $p/q$, the flow line $C$ has at most one
  outgoing end of any multiplicity at $e_{p/q}$.
\end{itemize}

Now suppose that $\phi$ is as in \S\ref{sec:cylinderResult} and that $J$ is
admissible and generic, but not necessarily local linear.  In this
case we define $\langle\delta\alpha,\beta\rangle$ to be the mod 2
count of flow lines $C\in\mc{M}(\alpha,\beta)/\R$ satisfying
conditions (i) and (ii) above.  Note that such $C$ are isolated by the
index calculation in \S\ref{sec:constraints}.

\begin{lemma}
\label{lem:IRC}
Let $\alpha$ and $\beta$ be generators of $C_*(X_1,X_2;P',Q')$ with no
periodic orbits in common, where $\alpha$ is obtained from $\beta$ by
rounding a corner.  Suppose $\phi$ is as in \S\ref{sec:cylinderResult}
with $Q\ge Q'$ and $J$ is admissible and generic.  Then:
\begin{itemize}
\item[\rm(a)]
$\langle\delta\alpha,\beta\rangle$ as above is well defined, i.e.\ the
set of $C\in\mc{M}(\alpha,\beta)/\R$ satisfying (i) and (ii) is finite.
\item[\rm(b)] $\langle\delta\alpha,\beta\rangle$ does not depend on
  $\phi$, $J$, or $Q$ as above.
\end{itemize}
\end{lemma}

\begin{proof}
  (a)\qua By Gromov compactness as in \cite[Lem.\ 9.8]{pfh2}, it is enough
  to show that there does not exist a $k$-times broken GFL from
  $\alpha$ to $\beta$ with $k\ge 1$ satisfying conditions (i) and (ii).
  
  Suppose that $C=(C_0,\ldots,C_k)$ is such a broken GFL.  Let $C_{i,j}$
  denote the components of $C_i$.  We have $I^{\op{vir}}(C_{i,j})=
  -\chi(C_{i,j})+\mu_\tau^0(C_{i,j})$, which together with the genus
  zero condition from (i) and equations \eqref{eqn:ellInd} and
  \eqref{eqn:hypInd} implies that
\begin{equation}
\label{eqn:3DV}
2e_-(C_{i,j})+h(C_{i,j})=2+I^{\op{vir}}(C_{i,j}).
\end{equation}
Also $\sum_{i,j}I^{\op{vir}}(C_{i,j})=1$.

Each $C_{i,j}$ must satisfy $I^{\op{vir}}(C_{i,j})\ge 0$.  Otherwise
\eqref{eqn:3DV} implies that $C_{i,j}$ has only one incoming end at
some $h_{p/q}$, and all outgoing ends elliptic.  In particular
the quasi-embedded curve underlying $C_{i,j}$ lives in a moduli space
of expected dimension $\le -1$ (even before modding out by the $\R$
action), which is impossible for generic $J$ (or even during a generic
one-parameter deformation).

Therefore one of the $C_{i,j}$'s has $I^{\op{vir}}=1$, while all other
$C_{i,j}$'s have $I^{\op{vir}}=0$.  Also, by
Proposition~\ref{prop:leftRight}, all but one of the $C_{i,j}$'s goes
between orbit sets with the same polygonal path, and hence by
Lemma~\ref{lem:slice} maps to a neighborhood of $\theta = p/q$ for
some $p/q$ depending on the $C_{i,j}$.  Let $\widehat{C}$ denote the
remaining $C_{i,j}$; this has two incoming ends corresponding to the
edges of the corner being rounded.

We claim that $I^{\op{vir}}(\widehat{C})=0$.  Suppose not, so that all
$C_{i,j}\neq\widehat{C}$ have $I^{\op{vir}}=0$.  Then
$\widehat{C}=C_k$.  (If $\widehat{C}=C_l$ with $l<k$, then each
$C_{i,j}$ with $i>l$ is a cylinder.  By \cite[Prop.\ 9.1]{pfh2}, any
cylinder with the same ends is trivial.  Thus one of these cylinders
must have distinct ends and hence $I^{\op{vir}}\neq 0$.)  Since $C$
has genus zero, downward induction on $i$ shows that each $C_{i,j}$
with $i<k$ has only one incoming end.  It follows by equation
\eqref{eqn:3DV} that each nontrivial $C_{i,j}$ with $i<k$ has more
than one outgoing end.  Since $k>1$, this leads to a contradiction of
condition (ii).

So $I^{\op{vir}}(\widehat{C})=0$, and since $\widehat{C}$ has incoming
ends of multiplicity one, $\widehat{C}$ is not multiply covered and
hence does not exist for generic $J$.

(b)\qua Consider a generic one-parameter deformation of $\phi$ and $J$.
By Gromov compactness as in \cite[Lem. 9.8]{pfh2},
$\langle\delta\alpha,\beta\rangle$ can change during the deformation
only at those times when there exists a broken GFL $C$ from $\alpha$
to $\beta$ satisfying (i) and (ii).  The classification of such broken
GFL's from part (a) is still valid, except that now $\widehat{C}$ as
above may exist at isolated times in a generic one-parameter family.
By equation \eqref{eqn:3DV}, the two incoming ends of $\widehat{C}$
are hyperbolic, while all incoming ends of $\widehat{C}$ are elliptic.
We claim that $k=1$, that $\widehat{C}=C_0$ or $\widehat{C}=C_1$, and
that all components of $C_1$ or $C_0$ respectively are trivial except
for one cylinder whose incoming end is at $e_{p/q}$ and whose outgoing
end is at $h_{p/q}$ for some $p/q$.  This follows from (i), (ii), and
\eqref{eqn:3DV} using induction over the components of $C$.  As in
Lemma~\ref{lem:MorseBott}, there are two cancelling cylinders from
$e_{p/q}$ to $h_{p/q}$.  Therefore standard gluing arguments as in
\cite{af,lee,sul} show that the mod 2 count
$\langle\delta\alpha,\beta\rangle$ does not change in this
bifurcation.

Finally, $\langle\delta\alpha,\beta\rangle$ is independent of $Q\ge
Q'$, because if $\phi$ satisfies the conditions in
\S\ref{sec:cylinderResult} for a given value of $Q$, then it also does
for any smaller value of $Q$.
\end{proof}

\subsection{$SL_2\Z$ symmetry}
\label{sec:symmetry}

We now observe a useful symmetry of our chain complex.  Let
$A=
\mbox{\begin{footnotesize}
$\begin{pmatrix}a&b\\c&d\end{pmatrix}$
\end{footnotesize}}
\in SL_2\Z$.  Define $I_A\eqdef \{x\in\R\mid cx+d>0\}$.
If $(p,q)$ is in the upper half plane and $p/q\in I_A$, then $A(p,q)$
is also in the upper half plane.  So if $[X_1,X_2]\subset I_A$, then
there is a well-defined linear map
\[
\Psi_A: CP_*(X_1,X_2;P,Q) \longrightarrow
CP_*\left(\frac{aX_1+b}{cX_1+d}, \frac{aX_2+b}{cX_2+d}; aP+bQ,cP+dQ\right)
\]
that replaces every occurrence of $e_{p/q}$ with
$e_{\frac{ap+bq}{cp+dq}}$ and $h_{p/q}$ with
$h_{\frac{ap+bq}{cp+dq}}$.  If $\alpha$ is obtained from
$\beta$ by rounding a corner, then $\Psi_A(\alpha)$ is obtained from
$\Psi_A(\beta)$ by rounding a corner.

\begin{lemma}
\label{lem:symmetry}
Under the assumptions of Lemma~\ref{lem:IRC}, let $A\in SL_2\Z$,
suppose $[X_1,X_2]\subset I_A$, and suppose $cP'+dQ'\le Q$.
Then
\[
\left\langle\delta\Psi_A(\alpha),\Psi_A(\beta)\right\rangle
=
\langle\delta\alpha,\beta\rangle.
\]
\end{lemma}

\begin{proof}
  Using the coordinates \eqref{eqn:cylinderCoordinates}, we define an
  orientation-preserving diffeomorphism of mapping tori
\[
\begin{split}
\psi_A:S^1\times[X_1,X_2]\times S^1
&\longrightarrow
S^1\times [cX_1+d,cX_2+d]\times
S^1,\\
(t,x,y) &\longmapsto
\left(-cy+dt,\frac{ax+b}{cx+d}, ay-bt\right).
\end{split}
\]
Then $\psi_A$ sends the mapping torus flow $R=\partial_t-x\partial_y$
for $\phi_0$ to a positive multiple of itself, namely
\[
(\psi_A)_*R = \frac{R}{a-cx}.
\]
In particular, $\psi_A$ sends the circle of periodic orbits at $x=p/q$
to the circle of periodic orbits at $x=(ap+bq)/(cp+dq)$.  Thus
$\psi_A$ pulls back a pair $(\phi',J')$ for $X_1,X_2,Q$ from an
admissible pair $(\phi,J)$ for $(aX_1+b)/(cX_1+d)$,
$(aX_2+b)/(cX_2+d)$, and $cP+dQ$.  But $J'$ is not admissible, because
$J'\partial_s$ is only a positive multiple of $R$, and $J'$ is
$\Psi_A^*\Omega$-tame but not necessarily $\Omega$-tame.

In the coordinates \eqref{eqn:cylinderCoordinates} we have
\[
\Omega = ds\wedge dt + dx\wedge dy + xdx\wedge dt,
\]
from which we compute that
\[
\Omega\wedge\Psi_A^*\Omega = \frac{(cx+d)^4+1}{(cx+d)^3}ds\wedge
dt\wedge dx\wedge dy > 0.
\]
Therefore linear interpolation defines a path of symplectic forms from
$\Omega$ to $\Psi_A^*\Omega$.  Let $J''$ be an admissible almost
complex structure for $X_1,X_2,Q$.  With respect to this path of
symplectic forms, we can find a path of tame almost complex structures
from $J''$ to $J'$ that are $\R$-invariant and that send the Reeb flow
to a positive multiple of itself.  Then a compactness argument as in
Lemma~\ref{lem:IRC} shows that the differential coefficient
$\langle\delta\alpha,\beta\rangle$ stays well defined and does not
change during the deformation, and hence equals
$\left\langle\delta\Psi_A(\alpha),\Psi_A(\beta)\right\rangle$.
\end{proof}

\subsection{Nonvanishing of some differential coefficients}
\label{sec:nonvanishing}

\begin{lemma}
\label{lem:nonvanishing}
Let $Q'\le Q$, and let $\alpha$ and $\beta$ be generators of
$C_*(X_1,X_2;P',Q')$ with no periodic orbits in common, where $\alpha$
is obtained from $\beta$ by rounding a corner.  Then for any $\phi$
and $J$ as in \S\ref{sec:cylinderResult}, we have
$\langle\delta\alpha,\beta\rangle = 1$.
\end{lemma}

\begin{proof}
More explicitly, we can write $\beta=e_{c/d}h_{a/b}$ or
$\beta=h_{c/d}e_{a/b}$ or $\beta=h_{c/d}h_{a/b}$, where
$(a,b)$, $(c,d)$ are lattice points in the upper half plane with
$\op{gcd}(a,b)=\op{gcd}(c,d)=1$ and $a/b<c/d$ and $b+d=Q'$.
Let
\[
E
\eqdef
E\left(\frac{a}{b}+\epsilon,\frac{c}{d}-\epsilon;a+c,b+d\right)
\]
where $0<\epsilon<1/Q^2$, and let $H$ be the sum of all generators
obtained by replacing an $e_{p/q}$ factor in $E$ with $h_{p/q}$.  In
the first two cases for $\beta$ we have $\alpha=E$, while in the third
case $\alpha$ is a summand of $H$.

To simplify notation, let $\delta^*$ be the dual differential defined
by
\[
\langle\alpha,\delta^*\beta\rangle \eqdef
\langle\delta\alpha,\beta\rangle.
\]
Then to prove the lemma it is
enough to show that
\begin{gather}
\label{eqn:keya}
\delta^*(e_{c/d}h_{a/b})=\delta^*(h_{c/d}e_{a/b})=E,\\
\label{eqn:keyb}
\delta^*(h_{c/d}h_{a/b})=H.
\end{gather}
We first observe that there can be no other terms in
$\delta^*(e_{c/d}h_{a/b})$, $\delta^*(h_{c/d}e_{a/b})$, and
$\delta^*(h_{c/d}h_{a/b})$ by Theorem~\ref{thm:rounding}(a).  Now write
\[
E=e_{p_1/q_1}^{\lambda_1}\cdots e_{p_k/q_k}^{\lambda_k}
\]
with $p_i/q_i>p_j/q_j$ for $i<j$.  We prove equations \eqref{eqn:keya}
and \eqref{eqn:keyb} in three steps.

\medskip{\bf Step 1}\qua
Suppose $k=1$.  Then we show that \eqref{eqn:keya} holds, and if also
$\lambda_1=1$, or equivalently $bc-ad=1$, then
\eqref{eqn:keyb} holds.

To do so, we reduce to the paper by Taubes \cite{t01} which studies
pseudoholomorphic thrice-punctured spheres on $\R\times S^1\times S^2$
for a certain almost complex structure.  By an appropriate identification
of $[X_1,X_2]\times S^1$ with a subset of $S^2$ contained between two
latitude lines, we deduce from
\cite[Thm.\ A.2]{t01} that there is an almost
complex structure $J_T$ on $\R\times Y_{\phi_0}\simeq \R\times
S^1\times [X_1,X_2]\times S^1$ such that:
\begin{itemize}
\item[\rm(i)]
$J_T$ is $\R\times S^1\times S^1$ invariant.
\item[\rm(ii)]
If $s$ denotes the
$\R$ coordinate, then $J_T$ sends $\partial_s$ to a positive multiple of
the mapping torus flow $R$.
\item[\rm(iii)]
$J_T$ is tamed by $\Omega$.
\item[\rm(iv)]
The moduli space $\mc{M}_T$ of $J_T$-holomorphic thrice-punctured
spheres with an outgoing end at $x=(a+c)/(b+d)$ and with incoming ends
at $x=a/b$ and $x=c/d$ consists of a single orbit of the $\R\times
S^1\times S^1$ action.
\end{itemize}

Let $\gamma_{p/q}$ denote the circle of periodic orbits of
$\phi_0$ at $x=p/q$.  Let
\[
\xi: {\mc{M}_T}/{\R} \longrightarrow
\gamma_{{p_1}/{q_1}}\times \gamma_{c/d} \times \gamma_{a/b} \simeq (S^1)^3
\]
denote the ``endpoint map'' sending a flow line to the periodic orbits
at its ends.  The $S^1$ actions on $\xi(\mc{M}_T/\R)$ by rotation in
the $-t$ and $y$ directions have weights $(p_1,c,a)$ and $(q_1,d,b)$
respectively.  It follows by (iv) that
\[
\pm[\xi(\mc{M}_T/R)] = (dp_1-cq_1)[\gamma_{p_1/q_1}\gamma_{c/d}] +
(bp_1-aq_1)[\gamma_{p_1/q_1}\gamma_{a/b}] +
(bc-ad)[\gamma_{c/d}\gamma_{a/b}]
\]
in $H_2(\gamma_{{p_1}/{q_1}}\times \gamma_{c/d} \times \gamma_{a/b})$.
Since the triangle with vertices $(0,0)$, $(c,d)$, and $(a+c,b+d)$
contains no lattice points other than the vertices and the
$\lambda_1-1$ lattice points on the edge in between $(0,0)$ and
$(a+c,b+d)$, we get
\begin{equation}
\label{eqn:endpointHomology}
\pm[\xi(\mc{M}_T/R)] 
= -[\gamma_{p_1/q_1}\gamma_{c/d}] +
[\gamma_{p_1/q_1}\gamma_{a/b}] +
\lambda_1[\gamma_{c/d}\gamma_{a/b}].
\end{equation}
We now perturb the setup so that $\gamma_{p/q}$ splits into an
elliptic orbit $e_{p/q}$ and a hyperbolic orbit $h_{p/q}$ for
$(p,q)=(a,b),(c,d),(p_1,q_1)$.  By Morse-Bott theory, cf.\ \cite{fb},
it follows from \eqref{eqn:endpointHomology} that there will be one
flow line from $e_{p_1/q_1}^{\lambda_1}$ with one outgoing end to each
of $h_{c/d}e_{a/b}$ and $e_{c/d}h_{a/b}$; and when $\lambda_1=1$,
there will be one flow line from $h_{(a+c)/(b+d)}$ to
$h_{c/d}h_{a/b}$.  (These flow lines are constructed by a gluing
argument, using the fact that $\M_T$ is cut out transversely by
\cite{t01,t04}.  There is only one flow line in each case by a
compactness argument as in in Appendix~\ref{app:MorseBott}.  These
flow lines are embedded by the adjunction formula \eqref{eqn:adj}.)

To complete the proof of Step 1, we can deform this perturbation of
$(\phi_0,J_T)$ to $(\phi,J)$ while preserving $\R$-invariance and
conditions (ii) and (iii) above.  As in Lemma~\ref{lem:IRC}, the mod 2
count of flow lines as above does not change during the deformation.

\medskip{\bf Step 2}\qua
We now prove equation \eqref{eqn:keya} when $k>1$, assuming it holds
whenever the positive integer $bc-ad$ is smaller, using Step 1 and
$(\delta^*)^2=0$ (which follows from $\delta^2=0$).

To set up the application of $(\delta^*)^2=0$, we need to introduce
several vectors in the $(p,q)$-plane.  We begin with
\[
(v,w)\eqdef\lambda_k(p_k,q_k)-(p_{k-1},q_{k-1}).
\]
We claim that $(v,w)$ points to the left of $(c,d)$, that is
$cw-dv>0$, i.e.
\[
\det\begin{pmatrix} c & \lambda_kp_k\\ d & \lambda_k q_k\end{pmatrix}
>
\det\begin{pmatrix} c & p_{k-1} \\ d & q_{k-1}\end{pmatrix}.
\]
Otherwise the lattice point
$(c,d)+(a,b)-2\lambda_k(p_k,q_k)$ would lie between $\mc{P}(E)$ and
$\mc{P}(e_{c/d}e_{a/b})$ or
(if equality holds and $k=2$ and $\lambda_1=1$)
on the line segment between $(0,0)$ and $(c,d)$, contradicting the
definition of $E$ or the assumption $\op{gcd}(c,d)=1$.  Similarly,
$(v,w)$ points to the left of $(a,b)$, as otherwise the lattice point
$(c,d)+(a,b)-(v,w)$ would lie between $\mc{P}(E)$ and
$\mc{P}(e_{c/d}e_{a/b})$ or on the line segment between $(c,d)$ and
$(c,d)+(a,b)$.  Hence there exists $A\in\op{SL}_2\Z$ sending
$(v,w)$ to the upper half plane while keeping $(a,b)$ and $(c,d)$ in
the upper half plane.  So by Lemma~\ref{lem:symmetry} we may assume
that $w>0$.  (By Lemma~\ref{lem:IRC}, we may assume that $Q$ is large
enough that Lemma~\ref{lem:symmetry} is applicable.)

Consider now the triangle with vertices $(c,d)$, $(c,d)+(a,b)$, and
$(c,d)+(a,b)-(v,w)$.  Choose a lattice point $(f,g)$ in this triangle
with minimal positive distance to the line through $(c,d)$ and
$(c,d)+(a,b)$.  Define
\[
\begin{split}
(a',b') & \eqdef (f,g) - (c,d),\\
(a'',b'') & \eqdef (a,b) - (a',b'),\\
(p',q') & \eqdef \lambda_k(p_k,q_k)-(a'',b'').
\end{split}
\]
The relevant vectors look something like this:

\begin{center}
{\small
\psfrag {ab}{$(a,b)$}
\psfrag {ab1}{$(a',b')$}
\psfrag {ab2}{$(a'',b'')$}
\psfrag {cd}{$(c,d)$}
\psfrag {fg}{$(f,g)$}
\psfrag {vw}{$(v,w)$}
\psfrag {cd}{$(c,d)$}
\psfrag {pq}{$(p',q')$}
\psfrag {pkqk}{$(p_{k-1},q_{k-1})$}
\psfrag {l}{$\lambda_k(p_{k},q_{k})$}
\includegraphics[width=7cm]{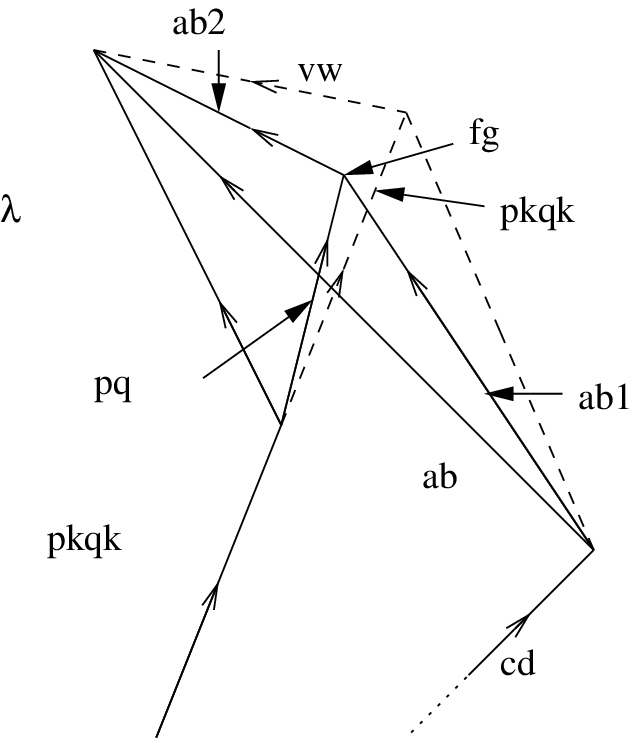}}
\end{center}

Note that $(p',q')$ points to the left of or in the same direction
as $(p_{k-1},q_{k-1})$, because the triangle with vertices $(c,d)$,
$(c,d)+(a,b)-(v,w)$, and
$(c,d)+(a,b)-\lambda_k(p_k,q_k)$ contains no lattice points other than
the vertices, because the same is true for
the triangle with vertices $(c,d)$, $(c,d)+(a,b)-\lambda_k(p_k,q_k)$,
and $(c,d)+(a,b)-\lambda_k(p_k,q_k)-(p_{k-1},q_{k-1})$, by definition of $E$.

Also $(p',q')$ points to the right of $(p_k,q_k)$, so $q'>0$.

The vector $(a',b')$ points to the right of $(a,b)$ and to the left of
$(c,d)$; in particular $b'>0$.

The vector $(a'',b'')$ points to the left of $(a,b)$ and to the right
of or in the same direction as $(v,w)$; hence $b''>0$.
Also, since $(a,b)$ and $(v,w)$ both point to the left of $(c,d)$, it
follows that $(a'',b'')$ does as well; hence
\begin{equation}
\label{eqn:inductionApplicable}
b'c-a'd<bc-ad.
\end{equation}
Now let $ E_0 \eqdef
\prod_{i=1}^{k-1}e_{p_i/q_i}^{\lambda_i} $.  Since there are no
lattice points between the paths $\mc{P}(E)$ and
$\mc{P}(e_{c/d}e_{a'/b'}e_{a''/b''})$, it follows that
\begin{gather*}
E\left(\frac{a'}{b'}+\epsilon,\frac{c}{d}-\epsilon;a'+c,b'+d\right)
=
E_0e_{p'/q'},\\
E\left(\frac{a''}{b''}+\epsilon,\frac{a'}{b'}-\epsilon;
a,b\right)
 = e_{a/b},\\
E\left(\frac{a''}{b''}+\epsilon,\frac{p'}{q'}-\epsilon;
a''+p',b''+q'\right)
= e_{p_k/q_k}^{\lambda_k}.
\end{gather*}
By Lemma~\ref{lem:inclusion}, we may assume that $X_1$ is small enough
so that $a''/b''\in[X_1,X_2]$. Then by inductive hypothesis using
\eqref{eqn:inductionApplicable}, and Step 1, we deduce
\begin{align}
\label{eqn:k-1}
\delta^*(e_{{c}/{d}}h_{{a'}/{b'}})&=
E_0e_{p'/q'},\\
\label{eqn:k=1a}
\delta^*(h_{a'/b'}h_{a''/b''})&=h_{a/b},\\
\label{eqn:k=1b}
\delta^*(e_{p'/q'}h_{a''/b''}) &=
e_{p_k/q_k}^{\lambda_k}.
\end{align}
By equations \eqref{eqn:k-1} and \eqref{eqn:k=1a}
and the trivial cylinder lemma~\ref{lem:trivCyl},
\begin{equation}
\label{eqn:delta*}
\delta^*(e_{c/d}h_{a'/b'}h_{a''/b''}) =
E_0e_{p'/q'}h_{a''/b''} + e_{c/d}h_{a/b}.
\end{equation}
No other terms are possible on the right hand side of
\eqref{eqn:delta*} by
Theorem~\ref{thm:rounding}(a).

It follows from \eqref{eqn:k=1b} and the trivial cylinder lemma
\ref{lem:trivCyl} that
\begin{equation}
\label{eqn:k=1c}
\left\langle E,\delta^*(E_0e_{p'/q'}h_{a''/b''})\right\rangle=1.
\end{equation}
Since $(\delta^*)^2=0$, equations \eqref{eqn:delta*} and
\eqref{eqn:k=1c} imply that
\[
\left\langle E,\delta^*(e_{c/d}h_{a/b})\right\rangle=1,
\]
so $\delta^*(e_{c/d}h_{a/b})=E$.

By a symmetric argument, $\delta^*(h_{c/d}e_{a/b})=E$.

\medskip{\bf Step 3}\qua
We now prove \eqref{eqn:keyb}.  Suppose first that $k=1$.  By Step 1
we may assume that $\lambda_1>1$.  By
\eqref{eqn:keya}, Theorem~\ref{thm:rounding}(a) and the trivial
cylinder lemma~\ref{lem:trivCyl} we have
\begin{equation}
\label{eqn:makesSense}
\delta^*(h_{c/d}e_{a/b}h_{a/b})
=e_{p_1/q_1}^{\lambda_1}h_{a/b}+x e_{p_1/q_1}^{\lambda_1-1}h_{p_1/q_1}e_{a/b},
\end{equation}
where $x$ is an unknown coefficient, and
\[
\delta^*\left(e_{p_1/q_1}^{\lambda_1}h_{a/b}\right)
= \delta^*\left( e_{p_1/q_1}^{\lambda_1-1}h_{p_1/q_1}e_{a/b} \right)
= e_{p_1/q_1}^{\lambda_1-1} e_{(a+p_1)/(b+q_1)}.
\]
By $(\delta^*)^2=0$ we get $x=1$, and by the trivial cylinder lemma
\ref{lem:trivCyl} and Theorem~\ref{thm:rounding}(a) we get
$\delta^*(h_{c/d}h_{a/b})=e_{p_1/q_1}^{\lambda_1-1}h_{p_1/q_1}$, so
\eqref{eqn:keyb} holds.

Strictly speaking, equation \eqref{eqn:makesSense} above makes sense
only if $Q$ is sufficiently large, e.g.\ if $Q\ge 2Q'$.  But we can
assume this without loss of generality by Lemma~\ref{lem:IRC}.

If $k>1$, we obtain \eqref{eqn:keyb} similarly to Step 2, using
$(\delta^*)^2(h_{c/d}h_{a'/b'}h_{a''/b''})=0$ and induction on
$bc-ad$.
\end{proof}

\begin{remark}
If we knew the expected isotopy invariance of PFH, cf.\
\S\ref{sec:pfh}, then we could give a conceptually simpler proof of
steps 1 and 2 above, without using \cite{t01}, by considering a Dehn
twist on a disc.  There we can easily calculate the PFH using isotopy
invariance, and together with a spectral sequence similar to the
one in \S\ref{sec:cylinderSpectral}, this gives sufficient
information about the differential.
\end{remark}

\begin{proof}[Proof of Theorem~\ref{thm:rounding}(b)]
This follows immediately from Lemma~\ref{lem:nonvanishing} and the trivial
cylinder lemma \ref{lem:trivCyl}.
\end{proof}

\subsection{Computing the homology of the cylinder complex}
\label{sec:cylinderSpectral}

\begin{proof}[Proof of Theorem~\ref{thm:cylinder}]
We compute the homology combinatorially using\break
Theorem~\ref{thm:rounding}.  We use induction. The base case of the
induction is when $P/Q\notin[X_1,X_2]$; here the theorem is obvious
because the chain complex has no generators.  For the inductive step,
suppose $P/Q\in[X_1,X_2]$.  We will prove the theorem assuming that it
is true for $(X_1,X_2';P',Q')$ whenever $X_2'<X_2$, $Q'\le Q$, and the
interval $[X_1,X_2']$ contains fewer rational numbers of denominator
$\le Q$ than the interval $[X_1,X_2]$.

To carry out the inductive step, define a filtration
$
\mc{F}_0\supset\mc{F}_{-1}\supset\cdots
$ of our chain complex as follows.  We say that the ``$p/q$ exponent'' of a
generator is the total exponent of $e_{p/q}$ and $h_{p/q}$.  Write
\begin{equation}
\label{eqn:recycle}
E=E(X_1,X_2;P,Q)=e_{p_1/q_1}^{\lambda_1}\cdots
e_{p_k/q_k}^{\lambda_k}
\end{equation}
with $p_1/q_1>\cdots>p_k/q_k$.  Let $c=p_1$ and $d=q_1$.  Define
$\mc{F}_{-i}$ to be the span of all generators with $c/d$ exponent at
least $i$.  We have
$
\delta(\mc{F}_{-i})\subset\mc{F}_{-i}
$
by Theorem~\ref{thm:rounding}(a), because if $p/q>c/d$, then
no generators have positive $p/q$ exponent; so rounding a
corner or double rounding never increases the $c/d$ exponent.

From this filtered complex we obtain a spectral sequence
$\mc{E}^*_{*,*}$ which converges to $HP_*(X_1,X_2;P,Q)$.  The
$\mc{E}^1$ term is the homology of the associated graded complex
$\mc{G}_{-i}=\mc{F}_{-i}/\mc{F}_{-i-1}$.  By
Theorem~\ref{thm:rounding}(a), if $\alpha$ and $\beta$ are generators
with the same $c/d$ exponent, then
$\langle\delta\alpha,\beta\rangle=1$ only if $\alpha$ is obtained from
$\beta$ by rounding a corner or double rounding not involving $c/d$, so
\[
\begin{split}
H_*(\mc{G}_0) &= HP_*(X_1,c/d-\epsilon;P,Q),\\ H_*(\mc{G}_{-i}) &=
\op{span}\left\{ e_{c/d}^i,e_{c/d}^{i-1}h_{c/d}\right\}
\tensor HP_*\left(X_1,\frac{c}{d}-\epsilon;P-ic,Q-id\right),\quad i>0.
\end{split}
\]
We now use the inductive hypothesis to make this more explicit.  Let
$l=\lambda_1$ above.  Then $H_*(\mc{G}_{-i})=0$ for $i>l$, while for
$i\le l$, the homology $
HP_*\left(X_1,{c}/{d}-\epsilon;P-ic,Q-id\right) $ is generated by the
homology class of
\[
E_i\eqdef E\left(X_1,\frac{c}{d}-\epsilon;P-ic,Q-id\right),
\]
together with the homology class $[H_i]$, where $H_i$ denotes a
generator obtained by replacing one of the $e_{p/q}$ factors in $E_i$
by $h_{p/q}$.

We now relate $E_i$ to $E_{i-1}$.
Suppose $0<i\le l$, and write
\[
E_i=e_{p_1/q_1}^{\lambda_1}\cdots e_{p_k/q_k}^{\lambda_k}\fedqe
e_{p_1/q_1}E_i'.
\]
with $p_1/q_1>\cdots>p_k/q_k$.  (Here we are recycling the notation so
that $p_j$, $q_j$, $\lambda_j$, $k$ are different from
\eqref{eqn:recycle}.)  Now $(p_1,q_1)$ is on the boundary of the convex hull
of the set of lattice points in the parallelogram $Z(X_1,c/d-\epsilon;
P-ic,Q-id)$.  It follows that $(p_1+c,q_1+d)$ is on the boundary of the
convex hull of the set of lattice points in
$Z(X_1,c/d-\epsilon;P-(i-1)c,Q-(i-1)d)$.  For example, the case $k=4$
and $\lambda_1=1$ looks like this:
\begin{center}
\begin{small}
\begin{picture}(300,250)(0,-30)

\put(60,0){\line(3,1){210}}
\put(60,0){\line(-1,2){60}}
\put(270,70){\line(-1,2){60}}
\put(0,120){\line(3,1){210}}

\put(210,50){\line(-1,2){60}}

\put(210,190){\line(2,-5){30}}
\put(240,115){\line(1,-6){5}}
\put(245,85){\line(-1,-1){25}}

\put(150,170){\line(2,-5){30}}
\put(180,95){\line(1,-6){5}}
\put(185,65){\line(-1,-1){25}}
\put(160,40){\line(-5,-2){100}}

\put(210,190){\circle*{3}}
\put(215,195){$\begin{pmatrix}P\\Q\end{pmatrix} -
(i-1)\begin{pmatrix}c\\d\end{pmatrix}$}

\put(150,170){\circle*{3}}
\put(90,190){
$\begin{pmatrix}P\\Q\end{pmatrix}-i \begin{pmatrix}c\\d\end{pmatrix}$}

\put(160,40){\circle*{3}}
\put(130,55){$\begin{pmatrix}p_1\\q_1\end{pmatrix}$}

\put(220,60){\circle*{3}}
\put(240,25){$\begin{pmatrix}p_1\\q_1\end{pmatrix} +
\begin{pmatrix}c\\d\end{pmatrix}$}

\put(240,40){\vector(-1,1){18}}

\put(120,20){\circle*{3}}
\put(122,1){$\begin{pmatrix}c\\d\end{pmatrix}$}

\put(60,0){\circle*{3}}
\put(62,-19){$\begin{pmatrix}0\\0\end{pmatrix}$}

\end{picture}
\end{small}
\end{center}
We deduce that
\[
E_{i-1} = E\left(\frac{p_1}{q_1}+\epsilon, \frac{c}{d}-\epsilon;
p_1+c, q_1+d\right)
E_i'.
\]
By Theorem~\ref{thm:rounding}(b), it follows that if we choose
$H_i=h_{p_1/q_1}E_i'$, then
\[
\begin{split}
\left\langle\delta\left(
e_{c/d}^{i-1}E_{i-1}
\right),
e_{c/d}^{i-1}h_{c/d}E_i
\right\rangle=
\left\langle\delta\left(
e_{c/d}^{i-1}E_{i-1}
\right),
e_{c/d}^iH_i
\right\rangle=1,\\
\left\langle\delta\left(
e_{c/d}^{i-2}h_{c/d}E_{i-1}
\right),
e_{c/d}^{i-1}h_{c/d} H_i
\right\rangle=1.
\end{split}
\]
Also, if we choose a representative $H_{i-1}$ of the class
$[H_{i-1}]$ by replacing one of the $e_{p/q}$ factors in
$E\left(p_1/q_1+\epsilon, c/d-\epsilon; p_1+c, q_1+d\right)$ by
$h_{p/q}$, then
\[
\left\langle\delta\left(
e_{c/d}^{i-1}H_{i-1}
\right),
e_{c/d}^{i-1}h_{c/d} H_i
\right\rangle=1.
\]
Finally, by Theorem~\ref{thm:rounding}(a), replacing $H_i$ with another
generator would change the corresponding differential coefficients
above to zero.

It follows that the first differential in the spectral
sequence, which we denote by $\partial_1$, satisfies
\begin{equation}
\label{eqn:cylinderDelta1}
\begin{split}
\partial_1\left(e_{c/d}^{i-1}[E_{i-1}]\right) &= e_{c/d}^{i-1}h_{c/d}[E_i] +
e_{c/d}^i[H_i],\\
\partial_1\left(e_{c/d}^{i-2}h_{c/d}[E_{i-1}]\right) &=
\partial_1\left(e_{c/d}^{i-1}[H_{i-1}]\right) = e_{c/d}^{i-1} h_{c/d} [H_i].
\end{split}
\end{equation}
The differential $\partial_1$ cannot contain any other terms due to the
bigrading, as we see by laying out the $\mc{E}^1$ term:
\begin{center}
\begin{picture}(330,280)(-270,-240)
\put(10,0){\framebox(60,40){$[H_0]$}}
\put(10,-40){\framebox(60,40){$[E_0]$}}

\put(-50,0){\framebox(60,40){$h_{c/d}[H_1]$}}
\put(-50,-40){\framebox(60,40){
$\begin{array}{c}e_{c/d}[H_1]\\h_{c/d}[E_1]\end{array}$}}
\put(-50,-80){\framebox(60,40){$e_{c/d}[E_1]$}}

\put(-115,-40){\framebox(65,40){$e_{c/d}h_{c/d}[H_2]$}}
\put(-115,-80){\framebox(65,40){
$\begin{array}{c}
e_{c/d}^2[H_2]\\
e_{c/d}h_{c/d}[E_2]
\end{array}$}}
\put(-115,-120){\framebox(65,40){$e_{c/d}^2[E_2]$}}

\put(-260,-150){\framebox(65,40){$e_{c/d}^{l-1}h_{c/d}[H_l]$}}
\put(-260,-190){\framebox(65,40){$\begin{array}{c}
e_{c/d}^l[H_l]\\
e_{c/d}^{l-1}h_{c/d}[E_l]
\end{array}$}}
\put(-260,-230){\framebox(65,40){$e_{c/d}^l[E_l]$}}

\put(30,20){\vector(-1,0){30}}
\put(30,-18){\vector(-4,1){30}}
\put(30,-22){\vector(-4,-1){30}}

\put(-37,-12){\vector(-4,-1){16}}
\put(-37,-28){\vector(-4,1){16}}
\put(-39,-58){\vector(-4,1){14}}
\put(-39,-62){\vector(-4,-1){14}}

\put(-104,-52){\vector(-4,-1){18}}
\put(-110,-68){\vector(-4,1){12}}
\put(-104,-98){\vector(-4,1){18}}
\put(-104,-102){\vector(-4,-1){18}}

\put(-150,-110){$\cdot$}
\put(-154,-113){$\cdot$}
\put(-158,-116){$\cdot$}

\put(-179,-122){\vector(-4,-1){18}}
\put(-179,-138){\vector(-4,1){18}}
\put(-180,-168){\vector(-4,1){18}}
\put(-180,-172){\vector(-4,-1){18}}

\end{picture}
\end{center}
From this picture we see that everything is killed in the spectral
sequence, except for
\[
e_{c/d}^l[E_l] = [E\left(X_1,X_2;P,Q\right)],
\]
and two generators which become homologous in $\mc{E}^2$, namely
\[
e_{c/d}^l[H_l], \;\; e_{c/d}^{l-1}h_{c/d}[E_l].
\]
Thus $HP_*(X_1,X_2;P,Q)$ is exactly as described in
Theorem~\ref{thm:cylinder}.
\end{proof}

\begin{remark}
The above algebraic calculation can be simplified after introducing some
more general combinatorial chain complexes involving rounding
corners of polygonal paths.  Compare \cite[Prop.\ 5.5]{t3}.
\end{remark}

\section{PFH of a Dehn twist on a torus}
\label{sec:torus}

Let $n$ be a positive integer.  We now consider the composition
$\phi_0^T$ of $n$ parallel positive Dehn twists on the torus from
equation \eqref{eqn:torusTwist}.  We can identify
\[
H_1\left(Y_{\phi_0^T}\right) \simeq \Z\oplus\Z\oplus\Z/n
\]
such that the circle of periodic orbits at $x=p/q$ is sent to
$(q,0,-p\op{mod} n)$.  Thus for each degree $d>0$, there are $n$
sectors $h\in H_1(Y)$ containing homology classes of orbit sets,
classified by the mod $n$ total numerator $[P]\in\Z/n$.  In this
section we compute $HP_*(\phi^T,h)$, where $\phi^T$ is a modification
of $\phi^T_0$ as in \S\ref{sec:cylinderResult} with $Q=d$, and $J$ is
any generic almost complex structure such that $(\phi^T,J)$ is
$d$-regular.  We denote this PFH by $HP_*(\phi^T;[P],d)$; this is
noncanonically $\Z/2d$-graded, see \S\ref{sec:pfh}.  Since the
isotopy from $\phi^T_0$ to $\phi^T$ is Hamiltonian, $(\phi^T,h)$ is
monotone, see Lemma~\ref{lem:monotone}.

\begin{theorem}
\label{thm:torus}
For every $[P]\in\Z/n$ and $d\in\Z^{>0}$, the periodic Floer homology
\[
HP_i\left(\phi^T;[P],d\right) \simeq \Z/2
\]
for each value of the $\Z/2d$-grading $i$.
\end{theorem}

In the proof we use the following notation.  As with the cylinder, for
$0<p/q<n$ and $q\le d$ there are periodic orbits $e_{p/q}$ and $h_{p/q}$
at $x=p/q$ of period $q$.  We denote the two fixed points at
$x=0\mod n$ simply by $e$ and $h$.  By symmetry, we may assume
without loss of generality that $P\equiv 0\mod n$.

\subsection{The wrapping spectral sequence}
\label{sec:eta}

The map $\phi^T_0$ has a circle of fixed points at $x\equiv 0 \mod n$,
each of which corresponds to a circle in $Y_{\phi^T_0}$.  The isotopy
from $\phi^T_0$ to $\phi^T$ induces a piecewise smooth homeomorphism
$Y_{\phi^T_0}\to Y_{\phi^T}$.  Let $\rho\subset Y_{\phi^T}$ be the
image of one of the circles in $Y_{\phi^T_0}$ coming from a fixed
point which does not survive under the perturbation.

\begin{definition}
If $C$ is a flow line for $(\phi^T,J)$, define the {\em wrapping number\/}
\[
\eta(C)
\eqdef
\#(C\cap(\R\times\rho))\in\Z.
\]
\end{definition}
Note that $\eta(C)$ depends only on the relative homology class of
$C$, and does not depend on $\rho$ since all such circles are
homotopic in the complement of the set of periodic orbits of period
$\le d$.  We remark that $\eta$ is reminiscent of the quantity $n_z$
considered in \cite{os2}.

Choose $0<\epsilon<1/d$ so that $\phi^T$ agrees with $\phi_0^T$ near
$x=\pm\epsilon\op{mod}n$.  If $\alpha$ is an orbit set, let
$d_0(\alpha)$ denote the total exponent of $e$ and $h$ in $\alpha$.

\begin{lemma}
\label{lem:eta}
If $C\in\M(\alpha,\beta)$, then $\eta(C)\ge 0$.
If $\eta(C)=0$, then:
\[
d_0(\beta)-d_0(\alpha) \ge
\left\{\begin{array}{cl} 1 & \mbox{if $C\cap\{x=\epsilon\}\neq\emptyset$}\\
0 & \mbox{otherwise}\end{array}\right.
+\;
\left\{\begin{array}{cl} 1 & \mbox{if $C\cap\{x=-\epsilon\}\neq\emptyset$}\\
0 & \mbox{otherwise.}\end{array}\right.
\]
\end{lemma}

\begin{proof}
Perturb $\epsilon$ so that $C$ is transverse to
$\{x=\pm\epsilon\}$.  Let
\[
(p_\pm,q_\pm)
\eqdef
[C\cap\{x=\pm\epsilon\}]
\in H_1(S^1_y\times S^1_t),
\]
with the sign conventions of \S\ref{sec:LEI}.  Evidently
\begin{align}
\label{eqn:peta}
p_\pm&=\eta(C),\\
\label{eqn:degreeFlux}
q_+-q_- &= d_0(\beta)-d_0(\alpha).
\end{align}
By Lemma~\ref{lem:slice}, we have $p_\pm \pm \epsilon q_{\pm}\ge 0$
with equality only if $C\cap\{x=\pm\epsilon\}=\emptyset$.  The lemma
follows.
\end{proof}

Since all flow lines have nonnegative wrapping number, we can write
\[
\delta=\delta_0+\delta_1+\delta_2+\cdots
\]
where $\delta_i$ counts the contributions from flow lines $C$ with
$\eta(C)=i$.  Since
$\delta^2=0$ and $\eta$ is additive under gluing of flow lines, it
follows that $\delta_0^2=0$, $\delta_0\delta_1=\delta_1\delta_0$, and
so forth.

By the index ambiguity formula \eqref{eqn:indexAmbiguity}, if
$C\in\M(\alpha,\beta)$, we can write
\begin{equation}
\label{eqn:torusIndex}
I(C)=I_0(\alpha,\beta)+2d\cdot\eta(C).
\end{equation}
In other words, our chain module has a relative $\Z$-grading $I_0$,
and $\delta_\eta$ shifts this grading by $2d\eta-1$.  We choose an
absolute $\Z$-grading by declaring the index of $e^d$ to be $0$.  In
particular, we have a $\Z$-graded $\delta_0$-homology $H_*(\delta_0)$,
and $\delta_1$ induces a map $(\delta_1)_*$ on it of degree $2d-1$.

By applying the following general algebraic lemma we obtain a spectral
sequence whose $\mc{E}^1$ term is $H_*(\delta_0)$, whose first
differential is $(\delta_1)_*$, and which converges to the PFH.  We
call this the {\em wrapping spectral sequence\/}.

\begin{lemma}
\label{lem:etaSS}
Let $C_*$ be a bounded $\Z$-graded vector space, and let
\[
\delta=\delta_0+\delta_1+\delta_2+\cdots:C_*\to C_*
\]
satisfy $\delta^2=0$ and $\op{deg}(\delta_i)=Ni-1$ with $N\neq 0$.
Then there is a spectral sequence $(\mc{E}^*,\hat{\delta}_*)$ such
that:
\begin{itemize}
\item $\mc{E}^r$ is $\Z$-graded and $\deg(\hat{\delta}_r)=Nr-1$.
\item
$\mc{E}^1_*=H_*(\delta_0)$, and
$\hat{\delta}_1=(\delta_1)_*:H_*(\delta_0)\to H_{*+N-1}(\delta_0)$.
\item
For $m\in\Z/N$ we have $H_m(\delta)\simeq \bigoplus_{i\equiv m
\op{mod} N}\mc{E}^\infty_i$.
\end{itemize}
\end{lemma}

\begin{proof}
  If $k$ is an integer, let $C_*[k]$ denote $C_*$ with the grading
  shifted by $k$.  Define a complex
  $\widetilde{C}_*\eqdef\bigoplus_{i\in\Z}C_*[Ni]$, with a filtration
\[
\mc{F}_i\widetilde{C}_*\eqdef \bigoplus_{j\le i}
C_*[Nj],
\]
where the component of the differential from $C_*[Ni]$ to $C_{*-1}[Nj]$ is
induced by $\delta_{i-j}$.  The filtration gives rise to a spectral
sequence $(\widetilde{\mc{E}}^*_{*,*},\widetilde{\delta}_*)$ with
$\widetilde{\mc{E}}^r_{p,q}\simeq \widetilde{\mc{E}}^r_{p-1,q-N+1}$.  We then
take $\mc{E}^r_q\eqdef \widetilde{\mc{E}}^r_{0,q}$, with $\hat{\delta}_r$
induced from $\widetilde{\delta}_r$.  Clearly
$(\mc{E}^*,\hat{\delta}_*)$ satisfies the first two properties;  since
$C_*$ is bounded and $N\neq 0$, this spectral sequence converges,
giving the third property.
\end{proof}

\subsection{Lifting from the torus to the cylinder}

To compute the differentials $\delta_\eta$, we need to relate flow
lines for a Dehn twist on a torus to flow lines for a Dehn twist on a
cylinder.  Under the covering
\[
\R\times S^1 \longrightarrow
(\R/n\Z)\times S^1
\]
of the torus by the infinite cylinder, the symplectomorphism $\phi^T$
of $(\R/n\Z)\times S^1$ lifts to a perturbation $\phi$ of the cylinder
twist $\phi_0$ on $\R\times S^1$.  There is then a covering of mapping
tori
\[
\pi:\R\times Y_\phi \longrightarrow \R\times Y_{\phi^T},
\]
and the almost complex structure $J$ for $\phi^T$ pulls back to an almost
complex structure $\pi^*J$ for $\phi$, such that $(\phi,\pi^*J)$
satisfy the conditions in \S\ref{sec:cylinderResult}.

The following lemma shows that to compute the PFH of $(\phi^T,J)$ on
the torus, we need only consider flow lines that lift to flow lines
for $(\phi,\pi^*J)$ on the cylinder.  (This lemma is actually used
only for a small part of the calculation in \S\ref{sec:degeneration}.)

\begin{lemma}
\label{lem:lifting}
Let $\delta'$ denote the contribution to the PFH differential for the
torus coming from flow lines that {\em do not\/} lift via $\pi$ to
flow lines for the cylinder.  Then $\delta'=0$.
\end{lemma}

\begin{proof}
Let $C$ be an $I=1$ flow line for the torus.  As in \S\ref{sec:index},
  since $I(C)=1$ we can write $C=C'\cup T$ where $C'$ is nontrivial
  and connected and $T$ is a union of trivial cylinders.
As in equation \eqref{eqn:pre3Dollar} we have
\begin{equation}
\label{eqn:TDT}
2g(C')+2e_-(C')+h(C')=3.
\end{equation}
In particular $g(C')\in\{0,1\}$.  Now $C$ lifts to the cylinder if and
only if $g(C')=1$.  For if $g(C')=0$, then $C'$ lifts to the cylinder
since each periodic orbit lifts.  Conversely, if $g(C')=1$, then
equation \eqref{eqn:TDT} implies that $C'$ is a flow line
from $\alpha$ to $\beta$, where $\alpha$ contains only elliptic
factors and $\beta=h_{p/q}$ for some $p,q$.  Such a $C'$ cannot lift
to the cylinder, since a lift would have $I<0$ by
Proposition~\ref{prop:area}, and hence does not exist for generic $J$.

By Lemma~\ref{lem:TIE}, for a certain (not locally linear) pair
$(\phi',J')$ where $\phi'$ is close to $\phi_0^T$, there are no genus
1 flow lines from $\alpha$ to $\beta$ as above.  Similarly to
Lemma~\ref{lem:IRC}, during a generic deformation of $(\phi',J')$ to
$(\phi^T,J)$, there are never any genus 1 broken GFL's from $\alpha$
to $\beta$ respecting the outgoing partitions at $\alpha$, so the mod
2 count of $I=1$ genus 1 flow lines from $\alpha$ to $\beta$ remains
zero.
\end{proof}

Going in the other direction, if $C$ is a flow line for the cylinder,
then it projects to a generalized flow line $\pi(C)$ for the torus.

\begin{lemma}
\label{lem:project}
If $C$ is a flow line for the cylinder and if $I(\pi(C))=1$, then
$\pi(C)$ is a flow line for the torus.
\end{lemma}

\begin{proof}
This follows immediately from Corollary~\ref{cor:restrict}.
\end{proof}

We now compute the wrapping number $\eta(\pi(C))$.  If
$\alpha=\gamma_1\cdots\gamma_k$ is an orbit set for the cylinder where
$\gamma_i=e_{p_i/q_i}$ or $\gamma_i=h_{p_i/q_i}$, define
\begin{equation}
\label{eqn:etaTilde}
\widetilde{\eta}(\alpha) \eqdef \sum_{i=1}^k
\left\lfloor\frac{p_i}{nq_i}\right\rfloor
\left(
-p_i + \frac{nq_i}{2}\left(
\left\lfloor\frac{p_i}{nq_i}\right\rfloor+1\right)\right).
\end{equation}

\begin{lemma}
\label{lem:projectionWrapping}
Let $C\in\mc{M}(\alpha,\beta)$ be a flow line for the cylinder.  Then
its projection to the torus has wrapping number
\[
\eta(\pi(C)) = \widetilde{\eta}(\alpha) - \widetilde{\eta}(\beta).
\]
\end{lemma}

\begin{proof}
  By equations \eqref{eqn:peta} and \eqref{eqn:sliceClass}, using the
  notation $p_i$, $q_i$, $p_j'$, $q_j'$, and $q_-$ from those
  equations, we have
\[
\begin{split}
(\eta(\pi(C)),q_-) &= [\pi(C)\cap\{x=-\epsilon\}]\\
&= \sum_{m\in\Z} \pi_*[C\cap\{x=mn-\epsilon\}]\\
&= \sum_{m\in\Z} \pi_*\left(\sum_{p_i/q_i>mn-\epsilon}(-p_i,q_i) -
\sum_{p_j'/q_j'>mn-\epsilon}(-p_j',q_j')\right)\\
&= \sum_{m\in\Z}\left( \sum_{p_i/q_i\ge mn}(-p_i+mnq_i,q_i) -
\sum_{p_j'/q_j'\ge mn}(-p_j'+mnq_j',q_j')\right).
\end{split}
\]
Note that the sum over $m$ has only finitely many nonzero terms.  The
lemma follows by a straightforward evaluation of this sum.
\end{proof}

The following lemma makes it easy to compute the differential $\delta_0$
in the wrapping spectral sequence.

\begin{lemma}
\label{lem:torusSplitting}
Lifting and projecting give a correspondence
\[
\left\{
\mbox{$I=1$, $\eta=0$ torus flow lines}\right\}
\longleftrightarrow
\left\{
\mbox{$I=1$ flow lines for $(-\epsilon,n+\epsilon)\times S^1$}
\right\}.
\]
This is a bijection, except that components of flow lines for the
torus in a neighborhood of $\{x=0\}$ can be lifted in two ways, to a
neighborhood of $\{x=0\}$ or to a neighborhood of $\{x=n\}$.
\end{lemma}

\begin{proof}
$(\rightarrow)$\qua Let $C$ be a flow line for the torus with $I(C)=1$ and
$\eta(C)=0$.  We need to show that $C$ lifts to an $I=1$ GFL for
$(-\epsilon,n+\epsilon)\times S^1$.  Clearly any GFL lifting $C$ is a
flow line, and unique modulo the proviso in the lemma.  Write
$C=C'\cup T$ where $C'$ is a nontrivial, connected flow line from
$\alpha$ to $\beta$, and $T$ is a union of trivial cylinders.

We first show that $C'$ lifts to a flow line for
$(-\epsilon,n+\epsilon)\times S^1$.  If $C'$ does not intersect both
of $\{x=\epsilon\op{mod} n\}$ and $\{x=-\epsilon\op{mod}n\}$, then
$C'$ trivially lifts to $(-\epsilon,n+\epsilon)\times S^1$.  If $C'$
intersects both these regions, then Lemma~\ref{lem:eta} implies that
$d_0(\beta)\ge 2$.  As in Lemma~\ref{lem:complexity}(c), any incoming
elliptic end of $C'$ has multiplicity one, so it follows from equation
\eqref{eqn:TDT} that $\beta=eh$ and $g(C')=0$, so $C'$ lifts to
$\R\times S^1$.  After translation we can choose the lift of $C'$ so
that $\beta$ is lifted to $h_{nk}e_0$ or $e_{nk}h_0$ with $k$ a
positive integer.  By Propositions~\ref{prop:area} and
\ref{prop:leftRight}, $\alpha$ is lifted to $e_{1/2}$ if $nk=1$, and
$e_{nk-1}e_1$ otherwise.  Lemma~\ref{lem:projectionWrapping} implies
that $\eta(C)=\eta(C')=k-1$.  Therefore $k=1$.  By
Lemma~\ref{lem:slice}, the lift is supported over
$(-\epsilon,n+\epsilon)\times S^1$.

We can also lift $T$ to obtain a lift of $C$.  The calculations in the
proof of Lemma~\ref{lem:QFormula} show that the lifted flow line
$\widetilde{C}$ satisfies $Q_\tau\big(\big[\widetilde{C}\big]\big)
=Q_\tau([C])$.  Together with equations \eqref{eqn:zeroChern},
\eqref{eqn:ellInd} and \eqref{eqn:hypInd}, this implies that
$I\big(\widetilde{C}\big)=I(C)$.

$(\leftarrow)$\qua Let $C$ be an $I=1$ flow line for
$(-\epsilon,n+\epsilon)\times S^1$; we need to show that the
projection $\pi(C)$ is a flow line for the torus with $I=1$ and
$\eta=0$. As above, $I(\pi(C))=I(C)$.  By
Lemma~\ref{lem:project}, $\pi(C)$ is a flow line.
Lemma~\ref{lem:projectionWrapping} implies that $\eta(\pi(C))=0$,
because if $\alpha$ is an orbit set for $(-\epsilon,n+\epsilon)\times
S^1$, then $\widetilde{\eta}(\alpha)=0$ by equation
\eqref{eqn:etaTilde} since $0\le p_i/nq_i\le 1$ for each $i$.
\end{proof}

\subsection{The $\eta=0$ homology}

We now compute the first term of the wrapping spectral sequence from
\S\ref{sec:eta}.

\begin{lemma}
\label{lem:torusSS}
For $d>0$ and $P\equiv 0 \op{mod} n$, the $\eta=0$ homology of a Dehn
twist on a torus is given by
\[
H_*(\delta_0) = \left\{\begin{array}{cl}
\Z/2 & \mbox{if $0\le * \le 2d-1$},\\
0 & \mbox{otherwise.}
\end{array}\right.
\]
\end{lemma}

\begin{proof}
We define a filtration on our chain complex by setting $\mc{F}_i$
to be the span of those generators in which the total exponent of $e$
and $h$ is at least $d-i$.  By Lemma~\ref{lem:eta} we have
$\delta_0(\mc{F}_i)\subset\mc{F}_i$.  Hence we obtain a spectral
sequence which converges to the $\eta=0$ homology.

The $\mc{E}^1$ term is given as follows.
Lemma~\ref{lem:torusSplitting} implies that the homology of the
associated graded complex is given in terms of the cylinder complex by
\[
\begin{split}
H_*(\mc{G}_0)&=\op{span}\{ e^d,e^{d-1}h\},\\
H_*(\mc{G}_i)&=\op{span}\{ e^{d-i},e^{d-i-1}h\}\tensor
\bigoplus_{k=1}^{i-1}HP_*(\epsilon,n-\epsilon;nk,i),\quad 2\le i\le d-1,\\
H_*(\mc{G}_d)&= \bigoplus_{k=1}^{d-1}HP_*(\epsilon,n-\epsilon;nk,d),
\end{split}
\]
and $H_*(\mc{G}_i)=0$ for all other $i$.
By Theorem~\ref{thm:cylinder}, if $0<k<i$, then the homology
$HP_*(\epsilon,n-\epsilon;nk,i)$ is generated by
the class
\[
E_{k,i}\eqdef 
[E(\epsilon,n-\epsilon;nk,i)],
\]
together with the homology class $H_{k,i}$ of those generators
obtained from the generator $E(\epsilon,n-\epsilon;nk,i)$ by replacing
an $e_{p/q}$ factor with $h_{p/q}$.

We now compute the first differential of the spectral sequence, which
we denote here by $\partial_1$.  By Lemma~\ref{lem:torusSplitting},
all contributions to $\partial_1$ come from projections to the torus
of differential coefficients in the cylinder complex
$CP_*(-\epsilon,n+\epsilon;*,d)$ which increase the total of the
$0$-exponent and the $n$-exponent by $1$.  When the $n$-exponent
increases by one, the differential coefficients are described by
\eqref{eqn:cylinderDelta1} with $X_1=\epsilon,X_2=n+\epsilon,c/d=n/1$.
When the $0$-exponent increases by one the differential coefficients have a
similar form by symmetry.  We conclude that
\begin{equation}
\label{eqn:yuck1}
\begin{split}
\partial_1(e^{d-i}E_{k,i})
&=
e^{d-i}h(E_{k,i-1}+E_{k-1,i-1})\\
&\;\;\;+e^{d-i+1}(H_{k,i-1}+H_{k-1,i-1}),\\
\partial_1(e^{d-i-1}hE_{k,i})& =\partial_1(e^{d-i}H_{k,i})\\
&=e^{d-i}h(H_{k,i-1}+H_{k-1,{i-1}}),\\
\partial_1(e^{d-i-1}hH_{k,i}) &=0.
\end{split}
\end{equation}
In this equation, when $k=1$ we interpret $E_{k-1,i-1}\eqdef
H_{k-1,i-1}\eqdef 0$, and when $k=i-1$ we interpret $E_{k,i-1}\eqdef
H_{k,i-1}\eqdef 0$.

For example, if
$n=2$ and $d=4$, then $(\mc{E}^1,\partial_1)$ is as follows:

\unitlength 0.985pt
\begin{picture}(380,195)(-270,-145)

\put(0,0){\framebox(110,50){$\begin{array}{c}
{[e_{3/2}h_{3/2}]} \\
{[e_{1/2}h_{3/2} = h_{1/2}e_{3/2}]} \\
{[e_{1/2}h_{1/2}]}
\end{array}$}}

\put(-110,0){\framebox(110,50){$\begin{array}{c}
h[e_1h_{3/2}=h_1e_{3/2}] \\
h[e_{1/2}h_1=h_{1/2}e_1] 
\end{array}$}}

\put(-170,0){\framebox(60,50){}}

\put(-220,0){\framebox(50,50){}}

\put(-270,0){\framebox(50,50){}}

\put(0,-65){\framebox(110,65){$\begin{array}{c}
{[e_{3/2}^2]} \\
{[e_{1/2}e_{3/2}]} \\
{[e_{1/2}^2]}
\end{array}$}}

\put(-110,-65){\framebox(110,65){$\begin{array}{c}
{e[e_1h_{3/2}=h_1e_{3/2}]} \\
{h[e_1e_{3/2}]} \\
{e[e_{1/2}h_1=h_{1/2}e_1]} \\
{h[e_{1/2}e_1]}
\end{array}$}}

\put(-170,-65){\framebox(60,65){$eh[e_1h_1]$}}

\put(-220,-65){\framebox(50,65){}}

\put(-270,-65){\framebox(50,65){}}

\put(0,-105){\framebox(110,40){}}

\put(-110,-105){\framebox(110,40){$\begin{array}{c}
e[e_1e_{3/2}]\\
e[e_{1/2}e_1]
\end{array}$}}

\put(-170,-105){\framebox(60,40){$\begin{array}{c}
e^2[e_1h_1]\\
eh[e_1^2]
\end{array}$}}

\put(-220,-105){\framebox(50,40){}}

\put(-270,-105){\framebox(50,40){$e^3h$}}

\put(0,-145){\framebox(110,40){}}

\put(-110,-145){\framebox(110,40){}}

\put(-170,-145){\framebox(60,40){$e^2[e_1^2]$}}

\put(-220,-145){\framebox(50,40){}}

\put(-270,-145){\framebox(50,40){$e^4$}}

\put(-270,-145){\framebox(380,195){}}

\put(29,41){\vector(-4,-1){37}}
\put(5,26){\vector(-4,1){13}}
\put(5,24){\vector(-4,-1){13}}
\put(29,9){\vector(-4,1){37}}

\put(38,-15){\line(-6,1){30}}
\put(8,-10){\vector(-1,0){18}}
\put(38,-17){\line(-6,-1){36}}
\put(2,-23){\vector(-1,0){32}}

\put(38,-48){\line(-6,1){30}}
\put(8,-43){\vector(-1,0){18}}
\put(38,-50){\line(-6,-1){36}}
\put(2,-56){\vector(-1,0){32}}

\put(32,-30){\line(-3,2){24}}
\put(8,-14){\vector(-1,0){18}}
\put(32,-32){\line(-6,1){30}}
\put(2,-27){\vector(-1,0){32}}
\put(33,-34){\line(-5,-1){25}}
\put(8,-39){\vector(-1,0){18}}
\put(32,-37){\line(-2,-1){30}}
\put(2,-52){\vector(-1,0){32}}

\put(-102,-10){\vector(-4,-3){20}}
\put(-102,-23){\vector(-4,-1){20}}
\put(-80,-23){\line(-1,0){22}}
\put(-102,-36){\vector(-4,1){20}}
\put(-98,-53){\vector(-4,3){24}}
\put(-80,-53){\line(-1,0){18}}

\put(-76,-75){\vector(-1,0){46}}
\put(-76,-78){\vector(-4,-1){48}}
\put(-76,-95){\vector(-1,0){48}}
\put(-76,-92){\vector(-4,1){46}}

\end{picture}

Here the notation $[a=b]$ means that one can choose $a$ or $b$, and
one will obtain the same homology class either way.

Equation \eqref{eqn:yuck1} implies that
the $\mc{E}^2$ term is given by
\begin{equation}
\label{eqn:yuck2}
\begin{split}
\mc{E}^2_{0,0}&=\op{span}\{ e^d\},\\
\mc{E}^2_{0,1}&=\op{span}\{ e^{d-1}h\},\\
\mc{E}^2_{i,i-2} &= \op{span}\left\{
e^{d-i}\sum_{k=1}^{i-1}E_{k,i} \right\}
,\quad 2\le i \le d,\\
\mc{E}^2_{i,i-1} &= \op{span}\left\{ 
\sum_{k=1}^{i-1}
e^{d-i-1}\left[hE_{k,i}=eH_{k,i}\right]
\right\},
\quad 2\le i \le d-1,\\
\mc{E}^2_{d,d-1} &= \op{span}\left\{
\sum_{k=1}^{d-1}H_{k,d}\right\},
\end{split}
\end{equation}
and all other $\mc{E}^2_{i,j}$'s are zero.  In equation
\eqref{eqn:yuck2}, strictly speaking there should be square brackets
around each generator to indicate that it is a homology class in
$\mc{E}^2$.  Here we have inferred the bigrading from the various
nonvanishing differential coefficients; this can also be computed
directly.

Equation \eqref{eqn:yuck2} shows that $\mc{E}^2$ has one generator in
each degree from $0$ to $2d-1$; so to complete the proof, we must show
that our spectral sequence degenerates at $\mc{E}^2$.  Because of the
bigrading, the only possible nonzero higher differential coefficient
is $ \left\langle \partial_2\left(e^{d-2}E_{1,2}\right), e^{d-1}h
\right\rangle $, which by the trivial cylinder lemma \ref{lem:trivCyl}
equals the $d=2$ differential coefficient
$\left\langle\partial_2\left(E_{1,2}\right),eh\right\rangle$.  But
this differential coefficient vanishes, because in the cylinder
complex $CP_*(-\epsilon,n+\epsilon;n,2)$ we have
$\delta(E_{1,2})=e_nh_0+h_ne_0$, which in the torus complex projects
to $eh+he=0$.
\end{proof}

\subsection{Degeneration of the wrapping spectral sequence}
\label{sec:degeneration}

We will now see that although there do exist index 1 flow lines with
$\eta>0$, these do not contribute in the wrapping spectral sequence.

\begin{proof}[Proof of Theorem~\ref{thm:torus}] We compute $H_*(\delta)$
using Lemma~\ref{lem:etaSS}.  Lemma~\ref{lem:torusSS} tells us that
$H_*(\delta_0)$ has one generator in each index $0,\ldots,2d-1$, and
we saw in the proof that
$[e^d]$ is the generator of index $0$.  Since $\delta_\eta$ shifts the
grading by $2d\eta-1$, it will suffice to show that
\[
(\delta_1)_*[e^d]=0 \in H_*(\delta_0).
\]
By Lemma~\ref{lem:lifting}, to compute $\delta_1$, it suffices to
consider flow lines that lift to the cylinder $\R\times {S^1}$.  As in
Lemma~\ref{lem:complexity}(c), the nontrivial component of $C$ has
only one outgoing end at $e$, so we can choose the lifted flow line
$C$ to be from $e_0^d$ to some $\beta$.

If $C$ is an $I=1$ flow line for the cylinder from $e_0^d$ to $\beta$,
then by Lemma~\ref{lem:project}, $\pi(C)$ contributes to
$\delta_1(e^d)$ if and only if $\eta(\pi(C)) = I(\pi(C))=1$.  By
equation \eqref{eqn:torusIndex} and
Lemma~\ref{lem:projectionWrapping}, this holds if and only if
$\widetilde{\eta}(\beta)=-1$ and $I_0(e^d,\pi(\beta))=1-2d$.
Therefore
\[
\delta_1(e^d) = \sum_{-\widetilde{\eta}(\beta)=I_0(e^d,\pi(\beta))+2d=1}
\langle\delta e_0^d,\beta\rangle \pi(\beta)
\]
where $\beta$ is a generator for the cylinder and $\pi(\beta)$ denotes
its projection to the torus.

In the cylinder, by Theorem~\ref{thm:rounding}, we have
\[
\delta\left(e_0^d\right)=
\sum_{k=1}^{d-1}\left(h_{\frac{1}{d-k}}e_{\frac{-1}{k}} +
e_{\frac{1}{d-k}}
h_{\frac{-1}{k}}\right)+\sum_{k=1}^{d-2}x_kh_{1/k}h_0h_{-1/(d-k-1)}
\]
where $x_1,\ldots,x_{d-2}\in\Z/2$ are unknown coefficients.  The
summands in the first sum on the right side all contribute to
$\delta_1(e^d)$.  To see this, we compute from equation
\eqref{eqn:etaTilde} that
\[
\widetilde{\eta}\left(h_{\frac{1}{d-k}}e_{\frac{-1}{k}}\right) =
-1.
\]
On the other hand, by
Lemma~\ref{lem:torusSplitting} and Proposition~\ref{prop:area},
\[
\begin{split}
I_0\left(e^d,h_{\frac{1}{d-k}}e_{\frac{nk-1}{k}}\right)
 & = 
I\left(e_0^{d-k}e_n^k,h_{\frac{1}{d-k}}e_{\frac{nk-1}{k}}\right)\\
&= 1-2d,
\end{split}
\]
where the right hand side denotes the index for the cylinder.
Likewise for $e_{\frac{1}{d-k}}h_{\frac{nk-1}{k}}$.  Thus
\[
\delta_1(e^d)=\sum_{k=1}^{d-1}\left(
h_{\frac{1}{d-k}}e_{\frac{nk-1}{k}}
+e_{\frac{1}{d-k}} h_{\frac{nk-1}{k}}\right) + h\left(\cdots\right).
\]
Unless $n=1$ and $k=1,d-1$, we have
\[
E(\epsilon,n-\epsilon;nk,d)=e_{\frac{1}{d-k}}e_{\frac{nk-1}{k}}.
\]
Then it follows by Theorem~\ref{thm:cylinder} that
$h_{1/(d-k)}e_{(nk-1)/k}$ and $e_{1/(d-k)} h_{(nk-1)/k}$ are
homologous in $CP_*(\epsilon,n-\epsilon;nk,d)$.  Hence, for all $n$
and $k$, $\delta_1(e^d)$ is $\delta_0$-homologous to a sum of
generators each containing an $e$ or an $h$.  So in terms of the
filtration in the proof of Lemma~\ref{lem:torusSS}, we have
$(\delta_1)_*[e^d]\in H_{2d-1}(\mc{F}_{d-1})$, but by
\eqref{eqn:yuck2}, $H_{2d-1}(\mc{F}_{d-1})=0$.
\end{proof}

\section{PFH of Dehn twists on higher genus surfaces}
\label{sec:surface}

Let $\Sigma$ be a compact connected symplectic surface, possibly with
boundary.  Choose a decomposition $\partial\Sigma=\partial_+\Sigma\sqcup
\partial_-\Sigma$.  Choose a finite number of {\em
disjoint\/} embedded circles $\gamma_i\subset\Sigma$, and to each
circle $\gamma_i$ associate a nonzero integer $n_i$.  In this
section we study the PFH in degree $d$ of the composition of $n_i$
positive Dehn twists along $\gamma_i$ for each $i$, for a perturbation
which is a small positive rotation on $\partial_+\Sigma$ and a small
negative rotation on $\partial_-\Sigma$.

\subsection{The setup}
\label{sec:surfaceSetup}

To be more precise, we define $\phi^\Sigma:\Sigma\to\Sigma$ as
follows.

First, let $N_i$ be disjoint tubular neighborhoods of the circles
$\gamma_i$ with coordinates $x_i\in[-\epsilon,|n_i|+\epsilon]$ and
$y_i\in\R/\Z$.  Here $0<\epsilon<1/d$.  On the cylinder $N_i$,
consider the twist
\[
\begin{split}
\phi_i^0:N_i
&\longrightarrow
N_i,\\
(x_i,y_i)
&\longmapsto
(x_i,y_i-x_i).
\end{split}
\]
If $n_i>0$, let $\phi_i$ be a perturbation of $\phi_i^0$ as in
\S\ref{sec:cylinderResult} with $Q=d$.
If $n_i<0$, let $\phi_i$ be the inverse of this perturbation.

Second, let $T_i\eqdef (\epsilon,|n_i|-\epsilon)\times\R/\Z\subset N_i$, and
let $\Sigma'\eqdef\Sigma\setminus\union_iT_i$.
Choose a Morse function
$
f:\Sigma'
\longrightarrow
[0,1]
$
such that
\[
\begin{split}
f^{-1}(1) &= \partial_+\Sigma\cup\union_{n_i>0}\partial \overline{T_i},\\
f^{-1}(0) &= \partial_-\Sigma\cup\union_{n_i<0}\partial \overline{T_i},
\end{split}
\]
and $|\nabla f|=1$ near the boundary.  Let $ \phi_f:\Sigma'
\longrightarrow \Sigma' $ be the time-1 Hamiltonian flow of $f$.
Choose a Riemannian metric on $\Sigma'$ which is large with respect to
$d$.  We assume that $f$ and the metric are Morse-Smale, and we let
$\partial^{\op{Morse}}$ denote the mod 2 differential in the Morse
complex.

We can make the above choices such that $\phi_i$ agrees with $\phi_f$
on $N_i\setminus T_i$, so the $\phi_i$'s and $\phi_f$ patch
together to give a symplectomorphism $\phi^\Sigma:\Sigma\to\Sigma$.

The period $\le d$ periodic orbits of $\phi^\Sigma$ consist of fixed
points at the critical points of $f$ in $\Sigma'$, as well as, in
each $T_i$, one elliptic and one hyperbolic orbit of period $q$ for
each rational number $p/q\in(0,|n_i|)$ with $q\le d$.

We choose a generic almost complex structure $J$ on $Y_{\phi^\Sigma}$,
such that $(\phi^\Sigma,J)$ is $d$-regular.  For convenience, we
assume that over $\Sigma'$, the almost complex structure $J$ is close
to the almost complex structure $J_0$ induced by the metric and
$\omega$ via the identification $Y_{\phi_f}\simeq S^1\times \Sigma'$
coming from the Hamiltonian isotopy.  (One can drop this assumption by
a modification of \S\ref{sec:splitting}.)

\begin{lemma}
\label{lem:monotone}
Under the assumption $(*)$ of \S\ref{sec:introduction}, if $\alpha$ is a
degree $d$ orbit set, then $(\phi^\Sigma,[\alpha])$ is monotone as in
\eqref{eqn:monotone}.
\end{lemma}

\begin{proof}
We have a short exact sequence
\[
0\longrightarrow H_2(\Sigma) \longrightarrow H_2\left(Y_{\phi^\Sigma}\right)
\longrightarrow \Ker\left(1-H_1\left(\phi^\Sigma\right)\right)
\longrightarrow 0.
\]
Assumption $(*)$ implies that
$\Ker\left(1-H_1\left(\phi^\Sigma\right)\right)$ is the image of the
inclusion-induced map $H_1(\Sigma')\to H_1(\Sigma)$.  The short exact
sequence then has a splitting
$\Ker\left(1-H_1\left(\phi^\Sigma\right)\right) \to
H_2(Y_{\phi^\Sigma})$, sending a loop $\xi\subset\Sigma'$ to
$S^1\times\xi\subset S^1\times\Sigma'\simeq Y_{\phi^f}$.  Since
$\phi^f$ is Hamiltonian isotopic to the identity on $\Sigma'$,
$[\omega]$ vanishes on the image of this splitting, as does
$c([\alpha])$.  So we just have to check that if
$\partial\Sigma=\emptyset$ then
\[
\int_{\Sigma}\omega = \lambda\langle[\Sigma],c([\alpha])\rangle.
\]
This holds if $d\neq g(\Sigma)-1$, because
$\langle[\Sigma],c([\alpha])\rangle = 2(d-g(\Sigma)+1)$.
\end{proof}

Thus under assumption $(*)$, we have a well defined
$\Z/2(d-g(\Sigma)+1)$-graded chain complex
\[
(CP_*(\phi^\Sigma,d),\delta) \eqdef \bigoplus_{h\cdot[\Sigma]=d}
(CP_*(\phi^\Sigma,h),\delta)
\]
whose differential $\delta$ may depend on $J$.

\subsection{The $\eta=(0,\ldots,0)$ complex}
\label{sec:splitting}

Continue to assume $(*)$.  We now describe a differential $\delta_0$ on
$CP_*(\phi^\Sigma,d)$, which is given explicitly in terms of Morse
theory on $\Sigma'$ and the cylinder complex for the $N_i$'s, and
which in some cases has the same homology as $\delta$.

Label the components of $\Sigma'$ by $\{\Sigma_j\mid
j=1,\ldots,\#\pi_0\Sigma'\}$.  For each $j$, let $\rho_j\subset
Y_{\phi^\Sigma}$ be a circle obtained from $S^1\times\{z_j\}\subset
S^1\times \Sigma_j$, where $z_j$ is not a fixed point of
$\phi^\Sigma$.  If $C$ is a flow line, define the {\em wrapping
number\/}
\[
\eta_j(C)\eqdef \#(C\cap (\R\times \rho_j))\in \Z.
\]
This does not depend on the choice of $z_j$.  Define
$\eta(C)\eqdef(\eta_1(C),\eta_2(C),\ldots)$.

We remark that if $\Sigma_j$ contains a component of $\partial\Sigma$,
then any flow line $C$ automatically has $\eta_j(C)=0$, since one can
choose $z_j$ near $\partial\Sigma$.

As in Lemma~\ref{lem:eta}, we have $\eta_j(C)\ge 0$.  Let $\delta_0$
denote the sum of the contributions to $\delta$ from flow lines $C$
with $\eta(C)=(0,\ldots,0)$.  Since $\eta_j(C)\ge 0$ and $\eta$ is
additive under gluing of flow lines, $\delta^2=0$ implies
$\delta_0^2=0$.

\begin{lemma}
\label{lem:splitting}
Let $\alpha$ and $\beta$ be generators of $CP_*(\phi^\Sigma,d)$.  Then
$\langle\delta_0\alpha,\beta\rangle=1$ if and only if either:
\begin{itemize}
\item[\rm(a)]
$\alpha=p\gamma$ and $\beta=q\gamma$, where
$p,q\in\op{Crit}(f)$
and
$\langle\partial^{\op{Morse}}p,q\rangle=1$, or:
\item[\rm(b)]
$\alpha=\alpha'\gamma$ and $\beta=\beta'\gamma$, where
$\alpha', \beta'$ are products of orbits in some $N_i$,
$\langle\delta\alpha',\beta'\rangle=1$ in $CP_*(\phi_i)$, and $\gamma$
is a product of orbits outside $N_i$.
\end{itemize}
\end{lemma}

\begin{proof}
We proceed in two steps.

\medskip{\bf Step 1}\qua Let $C$ be an index one flow line with
$\eta(C)=(0,\ldots,0)$ and without trivial cylinders.  We claim that $C$ is a
flow line for $\phi_f$ or for some $\phi_i$.

There is a canonical trivialization $\tau$ of $V$ over the $N_i$'s,
and also over the circle corresponding to each critical point of $f$.
There is no obstruction to extending this trivialization over
$\Sigma\setminus\{z_1,z_2,\ldots\}$, so $\eta(C)=(0,\ldots,0)$ implies
that $c_\tau([C])=0$.  It then follows as in \eqref{eqn:pre3Dollar}
that
\begin{equation}
\label{eqn:3DollarAgain}
2g(C)+2e_-(C)+h(C)=3.
\end{equation}
Here $h(C)$ denotes the total number of ends of $C$ at hyperbolic
orbits; and $e_-(C)$ denotes the number of incoming ends at elliptic
orbits in $T_i$'s with $n_i>0$ or at minima of $f$, plus the number of
outgoing ends at elliptic orbits in $T_i$'s with $n_i<0$ or at maxima
of $f$.

If our claim fails, then since $C$ is connected, WLOG there exists $i$
with $n_i>0$ such that $C$ intersects $N_i$ in both the regions
$\{x_i=-\epsilon/2\}$ and $\{x_i=+\epsilon/2\}$.  As in
Lemma~\ref{lem:eta}, $C$ must have incoming hyperbolic and elliptic
ends at $x_i=0$.  Then $C$ cannot intersect the region
$\{x_i=n_i-\epsilon/2\}$, or else as in Lemma~\ref{lem:eta} again, $C$
would have an incoming end at $x_i=n_i$, violating
\eqref{eqn:3DollarAgain}.  In the notation of equation
\eqref{eqn:sliceClass}, $[C\cap\{x=\epsilon/2\}]=(0,1)$, so by
equation \eqref{eqn:sliceClass}, $C$ has an outgoing end inside $T_i$.
Since this outgoing end is at $x_i>0$, equation \eqref{eqn:sliceClass}
again implies that there is another incoming end inside $T_i$,
contradicting \eqref{eqn:3DollarAgain}.

\medskip{\bf Step 2}\qua
We claim now that in $CP_*(\phi_f)$, if $\alpha$ and $\beta$ are
generators then $\langle\delta\alpha,\beta\rangle=1$ if and only if
(a) holds.  This follows from
standard arguments in Floer theory, cf.\ \cite{pfh1}.  Namely, if
$\gamma:\R\to\Sigma'$ is a gradient flow line of $f$, then
\[
S^1\times \op{graph}(\gamma)\subset S^1\times\R\times\Sigma'
\simeq \R\times Y_{\phi_f}
\]
is an embedded $J_0$-holomorphic cylinder cut out transversely, whose
PFH index agrees with the Morse index.  Index 1 flow lines not of this
form may exist for $(\phi_f,J_0)$, but $S^1$ acts nontrivially on
their moduli spaces, so by virtual cycle machinery as developed e.g.\
in \cite{cms,fo}, the mod 2 count of such bad flow lines is zero after
perturbing $J_0$ to our generic almost complex structure $J$.
\end{proof}

\subsection{A single nonseparating Dehn twist on a closed surface}
\label{sec:nonseparating}

Now let $\Sigma$ be a closed surface of genus $g$.  We consider a
single ($n=1$) positive Dehn twist along a nonseparating circle
$\gamma\subset\Sigma$, with neighborhoods $T\subset N$ and
$\Sigma'=\Sigma\setminus T$ as before.

\begin{theorem}
\label{thm:nonseparating}
If $\phi^\Sigma$ is a nonseparating positive Dehn twist on a closed
surface $\Sigma$ as above, and if $g\ge 2d+1$, then as
$\Z/2(d-g+1)$-graded modules,
\begin{equation}
\label{eqn:nonseparating}
HP_*(\phi^\Sigma,d) \simeq \Lambda^dK \oplus
\bigoplus_{d'=1}^d\Lambda^{d-d'}K\tensor
\{\mbox{$\Z/2$ in index $0,\ldots,2d'-1$}\},
\end{equation}
where $K\eqdef \op{Ker}(H_1(\Sigma',\partial\Sigma')\to
H_0(\partial\Sigma'))\tensor\Z/2$. 
\end{theorem}

\begin{proof}
We proceed in 2 steps.

\medskip{\bf Step 1}\qua
We first show that $H_*(\delta_0)$ agrees with the
right hand side of \eqref{eqn:nonseparating}.

Let $e_0,e_1,h_0,h_1$ denote the elliptic and
hyperbolic fixed points at $x=0,1$ in $N$.  Our Morse function $f$ on
$\Sigma'$ will have minima at $e_0$ and $e_1$ and saddle points at
$h_0$ and $h_1$.  We choose $f$ to have no other minima, a unique
maximum $m$ with $\partial^{\op{Morse}}(m)=h_0+h_1$, one saddle point
$s$ with $\partial^{\op{Morse}}(s)=e_0+e_1$, and $2(g-1)$
other saddle points in the kernel of $\partial^{\op{Morse}}$.

We define a filtration on the $\eta=0$ chain complex by setting
$\mc{F}_i$ to be the span of those generators containing at least $d-i$
factors corresponding to critical points of $f$ in $\Sigma'$.  This
gives rise to a spectral sequence with
\[
\mc{E}^1\simeq \bigoplus_{p=0}^d
HP_*\left(\phi_f,d-p\right)
\tensor HP_*\left(\phi^\Sigma|_T,p\right)
\]
by Lemma~\ref{lem:splitting}.  In the Morse complex of $f$,
$e_0$ and $e_1$ represent the same homology class $e$, and $h_0$ and $h_1$
represent the same homology class $h$, so
\[
HP_*\left(\phi_f,d-p\right)
\simeq
\Lambda^{d-p}K \oplus \bigoplus_{l=1}^{d-p}
\Lambda^{d-p-l}K\tensor \op{span}\left\{e^{l}, e^{l-1}h\right\}
\]
by Lemma~\ref{lem:splitting}.  It follows from the above two equations that
\[
\mc{E}^1\simeq \Lambda^dK\oplus\bigoplus_{d'=1}^d\Lambda^{d-d'}K
\tensor
\widehat{\mc{E}}^1(d'),
\]
where $\widehat{\mc{E}}(d')$ denotes the $\eta=0$ spectral sequence
for the torus from Lemma~\ref{lem:torusSS} with $n=1$, in degree $d'$.

By Lemma~\ref{lem:splitting}, the higher differentials on $\mc{E}$ are
determined by flow lines in $N$, and hence they are given by the
differentials on $\widehat{\mc{E}}$, tensored with the identity on the
$\Lambda^*K$ factors.  Together with Lemma~\ref{lem:torusSS}, this
proves our claim.

\medskip{\bf Step 2}\qua
The theorem is trivial if $d=0$, so assume $d>0$.  We now relate
$H_*(\delta)$ to $H_*(\delta_0)$ using the wrapping spectral sequence
from Lemma~\ref{lem:etaSS}.  If $C$ is a flow line from $\alpha$ to
$\beta$, we can write
\[
I(C)=I_0(\alpha,\beta)+ 2(d-g+1)\eta(C).
\]
Therefore $H_*(\delta_0)$ is $\Z$-graded, and if $\delta_i$ denotes
the contribution to $\delta$ from flow lines $C$ with $\eta(C)=i$,
then $\delta_i$ shifts the grading by $2(d-g+1)i-1$.  Thus
Lemma~\ref{lem:etaSS} is applicable as $d-g+1\neq 0$.  Since the
$\eta=0$ homology is supported in index $0,\ldots,2d-1$, the wrapping
spectral sequence will automatically degenerate at $\mc{E}^1$, giving
$H_*(\delta)\simeq H_*(\delta_0)$, provided that
\[
\left|2(d-g+1)i-1\right| > 2d-1
\]
for all $i>0$, i.e.\ if $g\ge 2d+1$.
\end{proof}

\subsection{A single separating Dehn twist on a closed surface}
\label{sec:separating}

Now let $\Sigma$ be a closed surface, and consider a single ($n=1$)
positive Dehn twist along a separating circle $\gamma\subset\Sigma$.
Departing slightly from the previous notation, denote the components
of $\Sigma'$ by $\Sigma_0$ and $\Sigma_1$, and let $g_j\eqdef g(\Sigma_j)$.

\begin{theorem}
\label{thm:separating}
If $\phi^\Sigma$ is a separating positive Dehn twist on a surface as
above, and if $g_0,g_1\ge 2d$, then as $\Z/2(d-g_0-g_1+1)$-graded modules,
\begin{equation}
\label{eqn:separating}
HP_*(\phi^\Sigma,d) \simeq \bigoplus_{q=0}^d\Lambda^{d-q}H_1(\Sigma')
\tensor
\bigoplus_{p=0}^q C_{p,q}[p^2-p(1-2g_0+2d)].
\end{equation}
Here, if $p\le q/2$, then $C_{p,q}$ has one generator in each degree
$0,2,\ldots,2p$, and also $2q-1,2q-3,\ldots,2(q-p)+1$ if $p>1$; and
$C_{p,q}=C_{q-p,q}$.
\end{theorem}

\begin{proof}
We proceed in three steps.

\medskip{\bf Step 1}\qua
We first compute the $\eta=(0,0)$ homology.

Our Morse function $f$ on $\Sigma$ will have saddle points at $h_0$
and $h_1$ and minima at $e_0$ and $e_1$.  We choose our labeling so
that $e_j,h_j\in\Sigma_j$.  We can assume that the only other
critical points of $f$ are a maximum $m_j\in\Sigma_j$ with
$\partial^{\op{Morse}}m_j=h_j$, and $2(g_0+g_1)$ saddle points with
$\delta^{\op{Morse}}=0$.

Note that if there exists an $\eta=(0,0)$ flow line from $\alpha$
to $\beta$, then for homological reasons, $\alpha$ and $\beta$ have
the same total numerator, where we regard critical points in
$\Sigma_j$ as having numerator $j$.  Let us redefine the saddle points
to have numerator zero; then since the saddle points are in the kernel
of $\delta^{\op{Morse}}$, it is still true by
Lemma~\ref{lem:splitting} that if
$\langle\delta_0\alpha,\beta\rangle=1$ then $\alpha$ and $\beta$ have
the same total numerator.  Thus we have a well-defined, relatively
$\Z$-graded subcomplex spanned by generators with total numerator $p$.
We denote its homology by $H_*(\delta_0)(p)$, and we choose an
absolute $\Z$-grading $I_0$ by declaring $I_0(e_0^{d-p}e_1^p)=0$.  We
claim then that, whether or not $g_0,g_1\ge 2d$, we have
\begin{equation}
\label{eqn:cpq}
H_*(\delta_0)(p) \simeq \bigoplus_{q=p}^d
\Lambda^{d-q}H_1(\Sigma')\tensor C_{p,q}.
\end{equation}
To prove \eqref{eqn:cpq}, it is
enough to do the case when $g_0=g_1=0$; the general case is obtained
by including the $2(g_0+g_1)$ saddle points with
$\partial^{\op{Morse}}=0$, which give rise to the
$\Lambda^*H_1(\Sigma')$ factors.

As in \S\ref{sec:nonseparating}, we compute $H_*(\delta_0)(p)$ using
the spectral sequence $\mc{E}(p)$ coming from the filtration by the
total degree in $\Sigma'$.  We have
\[
\mc{E}^1(p)=\bigoplus_{q=p}^d \bigoplus_{0\le r\le p,d-q}
HP_*(\epsilon,1-\epsilon;p-r,q)
\tensor \op{span}\{ e_0^{d-q-r}e_1^{r}\}.
\]
A somewhat lengthy calculation as in Lemma~\ref{lem:torusSS} then
establishes \eqref{eqn:cpq}.

\medskip{\bf Step 2}\qua
The theorem is trivial if $d=0$, so assume that $d>0$.  We claim now
that if $g_0,g_1\ge 2d$, and if $Z\in H_2(Y;\alpha,\beta)$ is a
relative homology class with $\eta(Z)=(\eta_0,\eta_1)\neq (0,0)$ and
$\eta_i\ge 0$, then
\[
I_0(\alpha)-I_0(\beta)\ge I(\alpha,\beta;Z)+2d-1.
\]
To prove the claim, WLOG we can take $\alpha=e_0^{d-p}e_1^p$ and
$\beta=e_0^{d-p-k}e_1^{p+k}$, with $0\le k\le d-p\le d$.  We must show
that $I(\alpha,\beta;Z)\le -2d+1$.    For the trivialization $\tau$ in
\S\ref{sec:splitting}, we have
\[
c_\tau(Z)=\eta_0(1-2g_0) + \eta_1(1-2g_1).
\]
A calculation as in equation \eqref{eqn:detFormula} shows that
\[
Q_\tau(Z)=k(d-2p-k) + d(\eta_0+\eta_1).
\]
Of course, $\mu_\tau(\alpha,\beta)=0$.  For homological reasons, we
have $\eta_0-\eta_1=k$.  Putting this all together, we get
\begin{equation}
\label{eqn:gradingShift}
I(\alpha,\beta;Z) = k(1-2g_0+2d-2p-k) + \eta_1(2-2g_0-2g_1+2d).
\end{equation}
If $k=0$, then $\eta_1>0$, so $g_0,g_1\ge 2d$ implies that
$I(\alpha,\beta;Z)\le -6d+2$.  If $k\ge 1$, then using $g_0\ge 2d$,
$p\ge 0$, $\eta_1\ge 0$, and $g_1\ge 2d>0$, we get
$I(\alpha,\beta;Z)\le -2d$.  Either way, our claim holds.

\medskip{\bf Step 3}\qua By Step 2, we
can write $\delta=\delta_0+\delta_1+\cdots$ where $\delta_i$ shifts
the $\Z$-grading $I_0$ on the $\eta=(0,0)$ complex by $-2i-1$, and
$\delta_1=\cdots=\delta_{d-1}=0$.  It follows from
Lemma~\ref{lem:etaSS} that the wrapping spectral sequence is well
defined with $N=-2$; and it automatically degenerates at $\mc{E}^1$, since the
$\eta=(0,0)$ homology is supported in degree $0,\ldots,2d-1$ by Step
1.  Thus, by \eqref{eqn:cpq}, we have
\begin{equation}
\label{eqn:z2Grading}
HP_*(\phi^\Sigma,d) \simeq \bigoplus_{q=0}^d\Lambda^{d-q}H_1(\Sigma')
\tensor
\bigoplus_{p=0}^q C_{p,q}
\end{equation}
as $\Z/2$-graded modules.  Now every differential in the wrapping
spectral sequence has degree $-1$ with respect to the
$\Z/2(d-g_0-g_1+1)$-grading.  By \eqref{eqn:gradingShift}, we can
refine \eqref{eqn:z2Grading} to an isomorphism of
$\Z/2(d-g_0-g_1+1)$-graded modules by inserting grading shifts as in
\eqref{eqn:separating}.
\end{proof}

\subsection{The general case}

For a composition of several Dehn twists along disjoint circles, as
considered in \S\ref{sec:surfaceSetup}, arguments like the above show
that the wrapping spectral sequence will exist and automatically
degenerate at $\mc{E}^1$, provided that condition $(**)$ from
\S\ref{sec:introduction} holds.  Then in principle the homology can be
computed combinatorially using Lemma~\ref{lem:splitting}, as
illustrated in the specific cases above.

Without condition $(**)$, one would need to understand the contribution
from flow lines with $\eta\neq (0,\ldots,0)$.  For all we know these
might never contribute to the differential,
except when some $\Sigma_j$ has genus $0$ and does not contain a
component of $\partial\Sigma$.

Without condition $(*)$, monotonicity often fails.  (Our composition of
Dehn twists will still have a monotone representative of its
symplectic isotopy class for each homology class $h$ as long as
$\partial\Sigma\neq\emptyset$ or $d\neq g-1$; however there might be
no such representative which is in standard form, i.e.\ equal to the
identity away from the twisting circles.) Without monotonicity one can
still define a version of PFH with coefficients in an appropriate
Novikov ring over $H_2(Y)$, which can be computed using the above
methods.

It is an interesting problem to attempt to extend these results to
compositions of Dehn twists along intersecting circles.  It would also
be interesting to try to compute PFH of pseudo-Anosov maps in terms of
hyperbolic geometry.

\begin{appendix}

\section{An argument from Morse-Bott theory}
\label{app:MorseBott}

Recall from \eqref{eqn:cylinderCoordinates} that we can identify
\[
\R\times Y_{\phi_0}\simeq \R\times S^1\times[X_1,X_2]\times S^1.
\]
We denote the coordinates on the right hand side by $s,t,x,y$.  There
is a natural almost complex structure $J_0$ on the right hand side
defined by
\begin{equation}
\label{eqn:j0}
\begin{split}
J_0(\partial_s)& \eqdef R,\\
J_0(\partial_x)& \eqdef \partial_y.
\end{split}
\end{equation}
Here $R= \partial_t-x\partial_y$ denotes the mapping torus flow.  Also
$J_0$ is invariant under the map
$(s,t,x,y)\mapsto(s,t,x+n,y-nt)$, and hence descends to $\R\times
Y_{\phi_0^T}$.

We can perform a perturbation of $(\phi_0,J_0)$ or $(\phi_0^T,J_0)$ to
obtain a pair $(\phi,J)$ such that $\phi$ satisfies the conditions in
\S\ref{sec:cylinderResult}, and $J$ is admissible (but $(\phi,J)$ is
not necessarily locally linear).  In this appendix, we use Morse-Bott
theory for the unperturbed setup to prove the following lemma.  This
explains Remark~\ref{rem:MorseBott} and is used in the proof of
Lemma~\ref{lem:lifting}.

\begin{lemma}
\label{lem:TIE}$\phantom{99}$
\begin{itemize}
\item[\rm(a)] There is a generic perturbation $(\phi,J)$ of
$(\phi_0,J_0)$ as above such that if $\alpha$ and $\beta$ are
generators of $CP_*(X_1,X_2;P,Q)$, and if $\alpha$ is obtained from
$\beta$ by double rounding, then genus $0$ flow lines in
$\mc{M}(\alpha,\beta)$ do not exist.
\item[\rm(b)] For every positive integer $q$, there is a generic
perturbation $(\phi,J)$ of $(\phi_0^T,J_0)$ as above, for which there
do not exist any genus $1$ flow lines with all outgoing ends elliptic
and with a single incoming end at a hyperbolic orbit $h_{p/q}$ of
multiplicity one.
\end{itemize}
\end{lemma}

\begin{proof}
(a)\qua For the unperturbed map $\phi_0$, let $\gamma_{p/q}$ denote the
circle of periodic orbits at $x=p/q$.  It will be convenient for
Lemma~\ref{lem:CQ} below to define an explicit diffeomorphism
$\theta:\gamma_{p/q}\to S^1$ by sending an orbit
$\gamma\in\gamma_{p/q}$ to
\begin{equation}
\label{eqn:identification}
\theta(\gamma)\eqdef qy+pt+\frac{pq}{2} \in \R/\Z
\end{equation}
for any point $(t,x=p/q,y)\in\gamma$.  This is well-defined since
$y+xt$ is constant on $\gamma$.

Suppose we are given rational numbers $p_1/q_1,\ldots,p_k/q_k$ (not
necessarily distinct), positive integers $\lambda_1,\ldots,\lambda_k$,
and also $p_1'/q_1',\ldots,p_l'/q_l'$ and
$\lambda_1',\ldots,\lambda_l'$.  Let $\widetilde{\mc{M}}$ denote the
moduli space of generalized flow lines for the unperturbed setup with
ordered ends, consisting of an outgoing end of multiplicity
$\lambda_i$ at some periodic orbit in the circle $\gamma_{p_i/q_i}$,
an incoming end of multiplicity $\lambda'_j$ at some periodic orbit in
the circle $\gamma_{p_j'/q_j'}$, and no other ends.  There is an
``endpoint map''
\[
\xi: \widetilde{\mc{M}} /\R \longrightarrow
\prod_{i=1}^k\gamma_{p_i/q_i} \times \prod_{j=1}^l\gamma_{p_j'/q_j'}
\]
which sends a flow line to the periodic orbits at its ends.  Given
$C\in\widetilde{\mc{M}}$, write
\[
\theta\circ\xi(C) =
(\theta_1,\ldots,\theta_k,\theta_1',\ldots,\theta_l') \in (S^1)^{k+l}.
\]

\begin{lemma}
\label{lem:CQ}
The endpoints of any $C\in\widetilde{\mc{M}}$ satisfy the linear
relation
\[
\sum_{i=1}^k\lambda_i\theta_i =
\sum_{j=1}^l\lambda'_j\theta'_j \in \R/\Z.
\]
\end{lemma}

\begin{proof}
It follows from \eqref{eqn:j0} that the 2-form $ dt\,dy - ds\,dx $ on
$\R\times Y_{\phi_0}$ annihilates any pair of tangent vectors of the
form $(v,J_0v)$.  Therefore
\[
\int_C dt\,dy = \int_Cds\,dx.
\]
Now $\int_Cds\,dx=0$ by Stokes theorem, because the 1-form $s\,dx$
vanishes along the periodic orbits.  Therefore $\int_Cdt\,dy=0$.  But
$\int_Cdt\,dy$ is just the area of the projection of $C$ to the
$(t,y)$-torus, and so under the identification
\eqref{eqn:identification},
\[
\int_Cdt\,dy \equiv
\sum_{i=1}^k\lambda_i\theta_i -
\sum_{j=1}^l\lambda'_j\theta'_j \mod \Z.
\]
The reason is that modulo $\Z$, the area of a homology in the
$(t,y)$-torus between a periodic orbit and a linear combination of the
curves $y=0$ and $t=0$ is equal to the right hand side of
\eqref{eqn:identification}.
\end{proof}

We now relate the perturbed setup to the unperturbed setup.  As in
\cite{fb}, we can perturb by choosing a Morse function $f_{p/q}$ on
$\gamma_{p/q}$ with one index $1$ critical point at $\tilde{h}_{p/q}$
and one index $0$ critical point at $\tilde{e}_{p/q}$.  Then
$\tilde{e}_{p/q}$ and $\tilde{h}_{p/q}$ get perturbed into $e_{p/q}$
and $h_{p/q}$.  Under the identification \eqref{eqn:identification},
it will be convenient below to choose
$\theta\left(\tilde{e}_{p/q}\right) = 0$ for all $p,q$, and set
$\theta\left(\tilde{h}_{p/q}\right)$ equal to some fixed irrational
number $\zeta$ between $2/3$ and $1$.

Suppose $(\phi_n,J_n)$ is a sequence of generic such perturbations
converging to the unperturbed setup, for which genus $0$ flow lines
$C_n$ exist from $\alpha$ to $\beta$.  We can pass to a subsequence
such that all the $C_n$'s have the same partitions at the ends.  Then
a compactness argument as in \cite{fb} shows that there exist GFL's
$C_0',\ldots,C_k'$ for the unperturbed setup such that:
\begin{itemize}
\item[\rm(i)]
There is a bijection between the outgoing ends of $C_0'$ and the
outgoing ends of each $C_n$.  If an outgoing end of each $C_n$ is at $e_{p/q}$,
then the corresponding end of $C_0'$ is at $\tilde{e}_{p/q}$.
\item[\rm(ii)]
For $0\le i < k$, there is a bijection between the incoming ends of
$C_i'$ and the outgoing ends of $C_{i+1}'$.  For each such pair, both
such ends are on the same $\gamma_{p/q}$.  Moreover, there is a
downward flow line of $f_{p/q}$ from the incoming end to the outgoing end.
\item[\rm(iii)]
There is a bijection between the incoming ends of $C_k'$ and the
incoming ends of each $C_n$.  If an incoming end of each $C_n$ is at
$h_{p/q}$, then the corresponding end of $C_k'$ is at
$\tilde{h}_{p/q}$.
\item[\rm(iv)]
The surface obtained by gluing the $C_i'$'s together along their
paired ends from (ii) has the same topological type as the $C_n$'s.
\end{itemize}

As in Proposition~\ref{prop:leftRight}, the polygonal path determined
by the outgoing ends of a GFL cannot cross to the right of the
polygonal path determined by the incoming ends. This allows us to
classify the possibilities for the curves $C_0',\ldots,C_k'$ as
follows.  Without loss of generality,
$\beta=h_{p_1'/q_1'}h_{p_2'/q_2'}h_{p_3'/q_3'}$, and $\alpha$ is a
product of $e_{p/q}$'s.  Furthermore the curves $C_0',\ldots,C_k'$
contain a total of at most two components that are not trivial
cylinders or branched covers thereof.  There cannot be just one such
component, because then that component would have all outgoing ends at
$\tilde{e}_{p/q}$'s and all incoming ends at $\tilde{h}_{p/q}$'s,
which contradicts Lemma~\ref{lem:CQ} since $\zeta$ is irrational.  So
there are two such components, call them $C_+$ and $C_-$.  We can
order the factors in $\beta$ so that $C_+$ has outgoing ends at some
$\tilde{e}_{p/q}$'s and incoming ends at $\tilde{h}_{p_1'/q_1'}$ and
some orbit $x_{p_0/q_0}\in\gamma_{p_0/q_0}$; while $C_-$ has outgoing
ends at some $\tilde{e}_{p/q}$'s and some orbit
$y_{p_0/q_0}\in\gamma_{p_0/q_0}$, and incoming ends at
$\tilde{h}_{p_2'/q_2'}$ and $\tilde{h}_{p_3'/q_3'}$.  By
Lemma~\ref{lem:CQ},
\begin{equation}
\label{eqn:thetaXi}
\begin{split}
\theta\left(x_{p_0/q_0}\right) &= -\zeta \in (0,1/3),\\
\theta\left(y_{p_0/q_0}\right) &= 2\zeta \in (1/3,\zeta).
\end{split}
\end{equation}
But by condition (ii) above, there is a downward flow line of
$f_{p_0/q_0}$ from $\theta\left(x_{p_0/q_0}\right)$ to
$\theta\left(y_{p_0/q_0}\right)$.  This contradicts
\eqref{eqn:thetaXi}, since $f_{p_0/q_0}$ takes its maximum at $\zeta$ and
its minimum at $0$.

(b)\qua For $\phi_0^T$, the map \eqref{eqn:identification} is not quite
well-defined, but it does give a well-defined two-to-one map
$\gamma_{p/q}\to \R/\frac{1}{2}\Z$.  Then Lemma~\ref{lem:CQ} still
holds in $\R/\frac{1}{2}\Z$.  In particular, if we choose
$\theta\left(\tilde{e}_{p/q}\right) = 0$ and set
$\theta\left(\tilde{h}_{p/q}\right)$ equal to some fixed irrational
number, then there do not exist any GFL's for $(\phi_0^T,J_0)$ with
all outgoing ends at $\tilde{e}_{p/q}$'s and all incoming ends at
$\tilde{h}_{p/q}$'s.

Suppose $(\phi_n,J_n)$ is a sequence of generic perturbations as above
converging to the unperturbed setup, for which genus $1$ flow lines
$C_n$ exist with a single incoming end at a hyperbolic orbit $h_{p/q}$
of multiplicity one, and all outgoing ends ellipitic.  Since $q$ is
assumed fixed, there are only finitely many possibilities for the ends
of $C_n$, so we can pass to a subsequence such that the $C_n$'s have
outgoing ends at the same elliptic orbits with the same
multiplicities, and an incoming end at the same hyperbolic orbit
$h_{p_-/q_-}$.  Then as before, there exist GFL's $C_0',\ldots,C_k'$
for the unperturbed setup satisfying conditions (i)--(iv) above.

Some component of some $C_i'$ must have genus $1$.  Otherwise all the
$C_i'$'s lift to the cylinder. As in Proposition~\ref{prop:leftRight},
if a GFL for the cylinder has one incoming end, then the outgoing ends
have the same underlying polygonal path.  So by condition (iii) and
downward induction on $i$, each $C_i'$ is a trivial cylinder on
$\tilde{h}_{p_-/q_-}$.  This leads to multiple contradictions.

Since the $C_n$'s also have genus $1$, it follows from condition (iv)
that each component of each $C_i'$ has only one incoming end.  The
genus zero components are then branched covers of trivial cylinders as
above.  This means that the genus $1$ component has an incoming end at
$\tilde{h}_{p_-/q_-}$ and all outgoing ends at $\tilde{e}_{p/q}$'s.
This contradicts our choices of $\theta\left(\tilde{e}_{p/q}\right)$
and $\theta\left(\tilde{h}_{p/q}\right)$.
\end{proof}

\end{appendix}

\Addresses\recd

\end{document}